\documentclass[reqno]{amsart}

\usepackage{xypic}
\input xy
\xyoption{all}
\usepackage{epsfig}
\usepackage{color}
\usepackage{amsthm}
\usepackage{amssymb}
\usepackage{amsmath}
\usepackage{amscd}

\newcommand{\labell}[1] {\label{#1}}

\numberwithin{equation}{section}
\newtheorem {Theorem}{Theorem} 
\numberwithin{Theorem}{section}

\newtheorem {Lemma}[Theorem]    {Lemma}         
\newtheorem {Proposition}[Theorem]{Proposition}  
\theoremstyle{definition}
\newtheorem{Definition}[Theorem]{Definition}
\theoremstyle{remark}
\newtheorem{Remark}[Theorem]{Remark}
\newtheorem{Example}[Theorem]{Example}

\newtheorem {Corollary}[Theorem]{Corollary}  

\expandafter\chardef\csname pre amssym.def at\endcsname=\the\catcode`\@
\catcode`\@=11
\def\undefine#1{\let#1\undefined}
\def\newsymbol#1#2#3#4#5{\let\next@\relax
 \ifnum#2=\@ne\let\next@\msafam@\else
 \ifnum#2=\tw@\let\next@\msbfam@\fi\fi
 \mathchardef#1="#3\next@#4#5}
\def\mathhexbox@#1#2#3{\relax
 \ifmmode\mathpalette{}{\m@th\mathchar"#1#2#3}%
 \else\leavevmode\hbox{$\m@th\mathchar"#1#2#3$}\fi}
\def\hexnumber@#1{\ifcase#1 0\or 1\or 2\or 3\or 4\or 5\or 6\or 7\or 8\or
 9\or A\or B\or C\or D\or E\or F\fi}

\font\teneufm=eufm10
\font\seveneufm=eufm7
\font\fiveeufm=eufm5
\newfam\eufmfam
\textfont\eufmfam=\teneufm
\scriptfont\eufmfam=\seveneufm
\scriptscriptfont\eufmfam=\fiveeufm

\catcode`\@=\csname pre amssym.def at\endcsname

\def    \eps    {\epsilon}

\newcommand{\CH}{{\mathcal H}}

\newcommand{\supp}{{\mathit supp}\,}
\newcommand{\point}{{\mathit point}}

\newcommand{\id}{{\mathit id}}

\newcommand{\Aa}{{\mathcal A}}

\newcommand{\Jj}{{\mathcal J}}

\newcommand{\Mm}{{\mathcal M}}

\newcommand{\Pp}{{\mathcal P}}

\def    \C      {{\mathbb C}}
\def    \R      {{\mathbb R}}
\def    \reals      {{\mathbb R}}
\def    \Z      {{\mathbb Z}}

\def    \CP     {{\mathbb C}{\mathbb P}}
\def    \ra     {{\rightarrow}}

\def    \12    {{\frac{1}{2}}}

\def    \codim  {\operatorname{codim}}

\def    \HF     {\operatorname{HF}}
\def    \CZ     {\operatorname{CZ}}
\def    \CF     {\operatorname{CF}}

\def    \ssminus        {\smallsetminus}

\def    \CW  {\operatorname{\bar{c}_{HZ}}}
\def    \wsh {\operatorname{\bar{c}_{hom}}}
\def    \CHZ  {\operatorname{c_{HZ}}}

\newcommand{\TCH}{\tilde{\mathcal H}}

\begin{document}


\setlength{\smallskipamount}{6pt}
\setlength{\medskipamount}{10pt}
\setlength{\bigskipamount}{16pt}





\title[Relative Hofer--Zehnder capacity]{Relative Hofer--Zehnder capacity
and periodic orbits in twisted cotangent bundles}

\author[Viktor Ginzburg]{Viktor L. Ginzburg}
\author[Ba\c sak G\"urel]{Ba\c sak Z. G\"urel}

\address{Department of Mathematics, UC Santa Cruz, 
Santa Cruz, CA 95064, USA}
\email{ginzburg@math.ucsc.edu, basak@math.ucsc.edu}

\date{\today}
\thanks{The work is partially supported by the NSF and by the faculty
research funds of the University of California, Santa Cruz.}

\bigskip

\begin{abstract}
The main theme of this paper is a relative version of the 
almost existence theorem for periodic orbits of autonomous Hamiltonian
systems. 

We show that almost all low levels of a function on a geometrically bounded 
symplectically aspherical manifold carry contractible periodic orbits of the 
Hamiltonian flow, provided  that the function attains its minimum
along a closed symplectic submanifold. As an immediate 
consequence, we obtain the existence of contractible periodic orbits on 
almost all low energy levels for twisted geodesic flows with symplectic 
magnetic field. We give examples of functions with a sequence of regular
levels without periodic orbits, converging to an isolated, but very 
degenerate, minimum.

The proof of the relative almost existence theorem 
hinges on the notion of the relative Hofer--Zehnder capacity and
on showing that this capacity of a small neighborhood of a symplectic 
submanifold is finite. The latter is carried out by proving that the flow of a 
Hamiltonian with sufficiently large variation has a non-trivial 
contractible one-periodic orbit, when the Hamiltonian is 
constant and equal to its maximum near a symplectic submanifold and 
supported in a neighborhood of the submanifold. 
\end{abstract}

\maketitle

\section{Introduction} 
In the framework of symplectic topology, Viterbo's proof of the Weinstein 
conjecture, \cite{vi1}, and its refinement, the almost existence theorem,
are among the most important results concerning the existence of periodic 
orbits of autonomous Hamiltonian systems. The almost existence theorem,
proved by H. Hofer and E. Zehnder, \cite{hz:book}, and by M. Struwe, 
\cite{St}, asserts that almost all (in the sense of measure theory) regular 
levels of a proper $C^2$-smooth function on $\R^{2n}$ carry periodic orbits 
of the Hamiltonian flow. In the last decade, these theorems
have been extended to a broad class of symplectic manifolds; see, e.g.,
\cite{fhv,hv-ct,hv,LT,lu,ma}. However, as has been
noted by E. Zehnder, the almost existence theorem fails for a general
symplectic manifold (see \cite{hz:book,Ze}).

In this paper, we prove a relative version of the almost existence theorem.
Namely, consider a function $F$ on a geometrically bounded symplectically
aspherical manifold $W$, with its minimum (say, equal to zero) attained along 
a closed symplectic submanifold $M$. The relative almost existence
theorem asserts that the levels $F=\eps$ carry 
contractible periodic orbits of the Hamiltonian flow of $F$ for almost all 
small $\eps>0$. This result strengthens a theorem of \cite{cgk} 
guaranteeing the 
existence of periodic orbits for a dense set of values of $F$ near zero.

This investigation has been inspired by the existence problem for periodic 
orbits of twisted geodesic flows (see Section \ref{subsec:magnetic} for 
definitions). As an immediate consequence of
the relative almost existence theorem, we obtain the existence of contractible
twisted geodesics on almost all low energy levels, provided that the magnetic
field is symplectic -- a result strengthening or
complementing numerous other results on the existence of twisted geodesics;
see, e.g., \cite{cgk,gi:survey,gk1,gk2,ke,lu,mac,pol2}. 

When $M$ is a Morse non-degenerate minimum of $F$ as is the case, for 
instance, for twisted geodesic flows, the flow 
has, hypothetically,  a periodic orbit on every low energy level. 
Moreover, one can expect
a certain lower bound on the number of such periodic orbits in terms of
the cohomology and codimension of $M$. This conjecture can be thought of
as a plausible generalization of the Weinstein--Moser theorem,
\cite{Mo:orbits,We:orbits}. (We refer the reader to \cite{gk1,ke} for
a detailed discussion of this conjecture and the proofs of 
certain particular cases.) The relative almost existence theorem provides
further evidence supporting the conjecture.

Note that unless the minimum of $F$ along $M$ is assumed to be
Morse--Bott non-degenerate, periodic orbits of the flow need not
exist on all levels of $F$ near the minimum. In  a variety of settings, 
we construct Hamiltonians $F$ such that for some sequence of regular values 
$\eps_k\to 0+$ the levels $F=\eps_k$ carry no
periodic orbits of the flow;  
see Section \ref{sec:counterexamples}. In these examples, 
zero is an isolated, but very degenerate, critical value. 
Such Hamiltonians also exist on the standard cotangent bundles to spheres.
These examples are obtained by combining A. Katok's construction, 
\cite{ka,zi}, the elimination of periodic orbits as in 
\cite{gi:bayarea,gg2,ke2}, and a suitable smoothing procedure, \cite{se};
see Sections \ref{sec:counterexamples} and \ref{sec:counterexamples-prfs}.

Similarly to the original almost existence theorem and to many other results
of this type, the proof of the relative almost existence theorem is comprised 
of two steps. 

The first step is showing that the flow of a Hamiltonian with sufficiently 
large variation has a non-trivial contractible one-periodic orbit. More 
precisely, we consider a Hamiltonian $H$ supported in a neighborhood of $M$, 
constant near $M$, and attaining its maximum along $M$. Then, if the maximum 
is large enough, the flow has a non-trivial contractible one-periodic orbit. 
This is proved by showing that the Floer homology of $H$ does not 
vanish for an interval of actions $[a,\,b)$ with $\max H < a$ 
(to ensure that the orbit is non-trivial). Note that establishing the 
existence of non-trivial one-periodic orbits is a common first step in
proving almost existence theorems. For instance, our theorem generalizes
a theorem proved in \cite{hz:book}. Similar theorems, under no
assumptions on $W$, but with some additional constraints on the normal
bundle to $M$ (e.g., that the normal bundle is trivial), have been
proved in \cite{hv,LT,lu2}. We also prove a (relatively easy) version of 
this theorem for time-dependent flows along the lines of \cite{bps}. Here, 
the orbits are still contractible, but not necessarily non-trivial.

The second step is to introduce the relative Hofer--Zehnder capacity. This
capacity is defined as in \cite{hz:book}, but for functions 
constant near $M$. The existence theorem for one-periodic orbits 
guarantees that the capacity of a small neighborhood of $M$ relative to $M$
is finite. The almost existence theorem follows now,
exactly as in \cite{hz:book}, from the fact that the capacity of the sets
$\{ F<\eps\}$ is an increasing function of $\eps>0$. 

The idea to consider the relative
Hofer--Zehnder capacity goes back to \cite{la}, and a number of existence
results akin to those proved here can be interpreted as calculations
of this capacity; see, e.g., \cite{hv,LT,lu2}.
Certain other versions of relative capacity have been recently introduced,
\cite{cgk,bps,lu2,sc}, but, to the best of our knowledge, the 
relative Hofer--Zehnder capacity has not been put to systematic use till now.
It should be noted that
it is not known whether or not this capacity differs from the original
Hofer--Zehnder capacity for neighborhoods of symplectic submanifolds
(see Example \ref{exam:CHZ-c0} and Section \ref{sec:other-cap}). Yet, 
in this setting, we can only establish finiteness of the relative capacity,
and hence this capacity readily lends itself as a convenient
tool for proving almost existence results.

The paper is organized as follows. In Section \ref{sec:main}, we state
the main results of the paper and define and briefly discuss the relative
Hofer--Zehnder capacity. The version of the theorem on the
existence of one-periodic orbits  for time-dependent Hamiltonians 
is stated in Section \ref{sec:time-dep}.
The rest of the paper is essentially devoted to the proofs.
The main goal of Section \ref{sec:prelim} is to recall some background
results and definitions needed for the proofs, in particular, those concerning 
filtered Floer homology. In Section \ref{sec:pfs}, we prove the main theorems
on the existence of one-periodic orbits (both the autonomous and time-dependent
cases). In Section \ref{sec:capacity-pf}, we  further discuss the relative 
Hofer--Zehnder capacity, compare it with some other capacities, and prove 
its properties  stated in Section \ref{sec:main}. The flows without
periodic orbits on a sequence of levels are constructed in Section 
\ref{sec:counterexamples-prfs}.

\subsection*{Acknowledgments.} The authors are deeply grateful to 
Paul Biran, Ely Kerman, Debra Lewis, Cesar Niche, Leonid Polterovich,
Felix Schlenk, and Ant\'ony Serra for useful discussions and
suggestions. The authors would also like to thank the Instituto
Superior T\'ecnico, Lisbon for its hospitality during the period when
a part of this work was carried out.

\section{Contractible periodic orbits of autonomous flows}
\labell{sec:main}
In this paper we are primarily concerned with almost existence theorems.
However, all known proofs of these theorems rely on first establishing
the existence of one-periodic orbits for Hamiltonians with sufficiently
large variation. Hence, we begin this section by stating some results
of this type.

\subsection{Periodic orbits of autonomous Hamiltonian systems}
The main theme of this subsection is the general principle asserting that 
under suitable additional hypotheses a compactly supported Hamiltonian
must have a \emph{fast non-trivial} periodic orbit, provided that the 
maximum of the Hamiltonian is sufficiently large. Here, we focus on 
Hamiltonians supported in a neighborhood of a closed symplectic submanifold 
$M$ and constrained on $M$. Let us start with a result 
which holds for a relatively broad class of Hamiltonians.

\begin{Theorem}
\labell{thm:min}
Assume that $M$ is a closed symplectic submanifold of a geometrically bounded
symplectically aspherical manifold $W$ (see the definitions below) and $U$
is a sufficiently small neighborhood of $M$. Then there exists
a finite constant $C>0$, depending on $U$, such that for every smooth
function $H$ supported in $U$ with $\min_M H> C$,
the Hamiltonian flow of $H$ has a non-trivial contractible periodic
orbit of period less than or equal to one.
\end{Theorem}

\begin{Remark}
When $W=M\times\R^{2n}$, equipped with the product symplectic structure,
the neighborhood $U$ can be taken arbitrarily large but bounded. For
$W=\R^{2n}$, Theorem \ref{thm:min} is due to H. Hofer and E. Zehnder;
see \cite[Theorem 12, p. 183]{hz:book}, where an explicit upper bound on
the period is given in terms of the capacity of $U$. 
A similar upper bound exists in terms of the relative capacity
of $(U,M)$; see Remark \ref{rmk:period}.

Note also that it is not known whether or not the function $H$ 
in Theorem \ref{thm:min} must have a contractible
orbit of period exactly equal to one.
\end{Remark}

Theorem \ref{thm:min}, 
proved in Section \ref{sec:capacity-pf}, is a soft consequence 
of a result which, under more restrictive
assumptions on $H$, allows one to control the actions on periodic orbits
and thus distinguish trivial and non-trivial orbits and establish the
existence of one-periodic orbits.
To state this result, we need first to recall some definitions
(including those used in Theorem \ref{thm:min}) and introduce
necessary notations.

For a compact subset $Z$ of a symplectic manifold $(V,\omega)$ without
boundary, denote by
$\CH(V,Z)$ the class of smooth functions $H\colon V\to \R$ such that
\begin{itemize}

\item[(H0)] $H$ vanishes on a non-empty open set (depending on $H$) whose
complement is compact;
\item[(H1)] $H$ is constant on a neighborhood of $Z$; and 

\item[(H2)] $\max_V H = H|_Z$.

\end{itemize}
When $V$ is not compact, the condition (H0) is equivalent to
requiring $H$ to have compact support. For a compact $V$, this condition
is equivalent to that $H$ vanishes on a neighborhood of a point.

If $V$ is a manifold with boundary and $Z$ is disjoint from
the boundary $\partial V$, we set $\CH(V,Z)=\CH(V\ssminus \partial V,Z)$. 
Note that then every function $H\in \CH(V,Z)$ extends to a smooth function 
on $V$ vanishing near $\partial V$.

Let $\omega|_{\pi_2(V)}=0$. Recall that for a time-dependent
Hamiltonian $H\colon S^1\times V\to \R$ the action functional 
on the space of smooth contractible loops is then defined as
\begin{equation}
\labell{eq:action}
A_H(x) = - \int_{D^2} \bar{x}^* \omega +
\int_{S^1}H(t,x)\,dt,
\end{equation}
where $x\colon S^1\to V$ is a contractible loop and $\bar{x}\colon D^2\to V$
is a map of a disk, bounded by $x$.

A smooth compactly supported function $H$ on an open subset $U$ of 
$V$ will always be regarded as a smooth compactly supported function on $V$ by 
extending it as zero to $V\ssminus U$. In particular, 
$\CH(U,Z)\subset C^\infty(V)$ when $Z\subset U$ and, for $H\in \CH(U,Z)$,
the action functional $A_H$ is defined on all contractible loops $x$ in $V$.

The main tool utilized in this paper to prove the existence of periodic 
orbits of Hamiltonian flows is Floer homology.
To have the Floer homology of $H$ defined, we need to impose some additional
conditions on the ambient manifold $V$ which will guarantee that the 
compactness theorem holds for $V$. The first of these conditions is that
$(V,\omega)$ is \emph{symplectically aspherical}, i.e.,
$$
\omega|_{\pi_2(V)}=0
\quad\text{and}\quad
c_1(TV)|_{\pi_2(V)}=0.
$$
Furthermore, since $V$ is not assumed to be compact, we need a way to
control the geometry of $V$ at infinity.
One standard way to do this is to require the manifold
to be convex at infinity. However, this requirement is not met in general 
by twisted cotangent bundles which serve as one of the motivating examples 
for this investigation. Hence, we impose a slightly weaker condition
and require the manifold to be geometrically bounded. Although the
precise definition is immaterial for what follows, we will recall it
for the sake of completeness.

\begin{Definition}
\labell{def:gb}
A symplectic manifold $(W,\omega)$ is said to be
\emph{geometrically bounded} if $W$ admits an almost complex structure $J$
and a complete Riemannian metric $g$ such that
\begin{itemize}
\item $J$ is uniformly $\omega$-tame, i.e.,
for some positive constants $c_1$ and $c_2$ we have
$$
\omega(X, JX)\geq c_1\| X\|^2 \quad\text{and}\quad
|\omega(X, Y)|\leq c_2\| X\|\,\|Y\|
$$
for all tangent vectors $X$ and $Y$ to $W$.

\item The sectional curvature of $(W,g)$ is bounded from above and
the injectivity radius of $(W,g)$ is bounded away from zero.
\end{itemize}
\end{Definition}

We refer the reader to
Chapters V (by J.-C. Sikorav) and X (by M. Audin, F. Lalonde, and L. 
Polterovich) in \cite{al} and to \cite{cgk,lu} for a discussion of the 
concept of geometrically bounded manifolds. 
Here we only mention that among these manifolds are compact manifolds,
manifolds convex at infinity (e.g., $\C^n$), and twisted cotangent bundles.

In what follows, we will always denote a geometrically bounded symplectically
aspherical manifold by $W$, while $V$ will stand for a general symplectic 
manifold. Likewise, $M$ will denote a closed (symplectic, in many instances)
submanifold of $W$ or $V$ and
$Z$ will be just a compact subset. 

Now we are in a position to state the main technical result of this section
which is the key to the proof of Theorem \ref{thm:min} and to the
relative almost existence theorem discussed later.

\begin{Theorem}
\labell{thm:main}
Assume that $M$ is a closed symplectic submanifold of a geometrically bounded
symplectically aspherical manifold $W$ and $U$ is
a sufficiently small neighborhood of $M$. Then there exists
a finite constant $C>0$, depending on $U$, such that for every $H\in \CH(U,M)$
with $\max H> C$ the Hamiltonian flow of $H$ has a non-trivial contractible
one-periodic orbit with action greater than $\max H$.
\end{Theorem}

This theorem will be proved in Section \ref{sec:pfs}.

\begin{Remark}
\labell{rmk:$C$}
When $W=M\times\R^{2n}$, equipped with the product symplectic structure,
in Theorem \ref{thm:main}, as in Theorem \ref{thm:min},
the neighborhood $U$ can be taken arbitrarily large but bounded.
The proof of Theorem \ref{thm:main} gives also an upper bound on the value
of the constant $C$ (see Theorem \ref{thm:cap-nbd}). Namely, it is 
sufficient to take $C=\pi R^2$, where
$R$ is the radius of a symplectic tubular neighborhood of $M$ containing $U$;
see Section \ref{sec:tubular} for the definition. Note also that the
constant $C$ in Theorem \ref{thm:min} can be taken the same as in Theorem
\ref{thm:main}. Furthermore, the function $H$ from
Theorem \ref{thm:main} must have non-trivial periodic orbits for every
period $T\geq 1$. (Indeed, the orbits of period $T\geq 1$ for $H$ are exactly
one-periodic orbits for the function $T\cdot H$ which clearly satisfies
the hypotheses of the theorem.)
\end{Remark}

\begin{Example}
Let $W=\R^{2n}$ be equipped with the standard symplectic structure. Then $M$ 
is necessarily a point. In this case, Theorem \ref{thm:main} is established 
by H. Hofer and E. Zehnder in \cite{hz:book} under the additional assumption 
that $H$ is non-negative.

When the normal bundle to $M$ in $W$ is trivial, a result similar to
Theorem \ref{thm:main} was proved by H. Hofer and C. Viterbo by a different 
method, \cite{hv}. In \cite{hv}, the manifold $W$ need not be symplectically 
aspherical, but the function $H$ is required to be non-negative.
\end{Example}

\begin{Remark}
In Theorem \ref{thm:main} and in the almost existence theorems stated below, 
the assumptions that $W$ is geometrically bounded and symplectically 
aspherical appear to be superfluous.
When $H\geq 0$, it should be possible to prove a version of the theorem 
without these assumptions, e.g., by utilizing the methods of \cite{hv,LT,lu}.
(Some preliminary results in this direction have been obtained by
E. Kerman, \cite{Ke:new}, and L. Macarini,   \cite{mac2}.)
However, we feel that the Floer homology proof given here is of interest 
by itself (even for $W=\R^{2n}$, at least because it allows one to eliminate
the assumption $H\geq 0$) and this argument may have some other applications.
\end{Remark}

As stated, Theorem \ref{thm:main} does not hold when the assumption (H1) 
that $H$ is constant near $M$ is replaced by the less restrictive assumption 
that $H$ is constant on $M$; see Example \ref{exam:actions}.
However, the hypotheses of Theorem \ref{thm:main} can be replaced
by the weaker condition $\min_M H >C$, as Theorem \ref{thm:min} indicates,
at the cost of loosing control of the value of the action (and, to some
extent, of the period of the orbit). Then, the flow still has a non-trivial 
contractible periodic orbit with period not exceeding one, but the action on 
this orbit does not have to be greater than or equal to $\max H$.

Furthermore, as will be clear from the proof, the assumption (H1) can be 
replaced by that all partial derivatives of $H$ at $M$ vanish up to fourth 
order. One can also allow $H$ to be time-dependent. In this case, we need to
require $H_t$, $t\in [0,1]$, to be periodic in time and belong to the class 
$\CH(U,M)$ for every $t$, and the maximum of $H$ has to be independent of 
time. It is not clear whether or not this last condition can be relaxed. 

Finally note that Theorem \ref{thm:main} does not extend to arbitrarily 
large bounded neighborhoods of $M$ if periodic orbits are still required 
to be contractible; see Example \ref{exam:horocycle}.

\subsection{Relative Hofer--Zehnder capacity and almost existence}
\labell{sec:HZ}

A relative version of the Hofer--Zehnder capacity can be defined in
variety of ways depending on whether or not (and how) the homotopy class
of an orbit is incorporated into the definition. 

\subsubsection{Relative Hofer--Zehnder capacity}
\labell{sec:defHZ}
The simplest definition of the capacity imposes no requirement on the
homotopy class of periodic orbits.

Let $(V,\omega)$ be a symplectic manifold.
A non-trivial periodic orbit of a Hamiltonian flow on $V$ will
be called \emph{fast} if the orbit has period less than or equal to one. When
the period is greater than one we will call the periodic orbit \emph{slow}.
The next definition, in which we follow \cite{la}, is the key to deriving
the almost existence theorem from Theorem \ref{thm:main}.

\begin{Definition}
\labell{def:CHZ}
Let $Z$ be an arbitrary compact non-empty subset of $V$. 
The \emph{relative Hofer--Zehnder capacity}
$\CHZ(V,Z)\in (0,\infty]$ is defined as
$$
\CHZ(V,Z)=\sup_{ H\in \CH(V,Z)}\{\max H\mid 
\text{all non-trivial periodic orbits of $H$ are slow}\}.
$$
When the dependence of this capacity on $\omega$ is essential, we will use
the notation $\CHZ(V,Z,\omega)$.
\end{Definition}

Thus, a function $H\in \CH(V,Z)$ with $\max H> \CHZ(V,Z)$ must have a 
non-trivial  fast periodic orbit.

As we have pointed out above, the assumptions (H1) and (H2) in the definition 
of the class $\CH(V,Z)$ can be replaced by the weaker conditions. Namely,
for $C>0$, set
$$
\TCH_C(V,Z)=\{ H\in C^\infty(V)\mid \text{$H$ satisfies (H0) and
$\min_Z H > C$}\}.
$$ 

\begin{Theorem}
\labell{thm:HZ-modify}
~\begin{enumerate}
\item $\CHZ(V,Z)=\inf C$, where the infimum is taken over $C>0$ such that
every $H\in \TCH_C(V,Z)$ has a non-trivial fast periodic orbit.

\item The capacity $\CHZ(V,Z)$ does not change when functions in $\CH(V,Z)$ or
$\TCH_C(V,Z)$ are required to be non-negative.
\end{enumerate}

\end{Theorem}

The proof of this theorem, nearly identical to the proof of
Theorem \ref{thm:min}, will be given in Section \ref{sec:capacity-pf}.

\begin{Corollary}
$\CHZ(V,\point)$ is equal to $\CHZ(V)$, the standard Hofer--Zehnder capacity 
of $V$ (see \cite{hz:book} for the definition).
\end{Corollary}

The main point of this corollary is that the conditions on the function $H$
in the definition of the standard Hofer--Zehnder capacity can be relaxed,
\cite[p.\ 70]{hz:book}.
Namely, it is not necessary to assume that $H$ vanishes on an open set
(rather at one point) or that $H$ attains its
maximal value on the complement to a compact set (only that $H$ is constant
outside a compact set). Strangely, this fact is not mentioned in \cite{hz:book}
although all the ingredients needed for the proof are already there.

The properties of the relative capacity are summarized in the following

\begin{Theorem}
\labell{thm:CHZ}
~
\begin{enumerate}

\item {\rm [Invariance].} The relative capacity $\CHZ$ is an invariant of 
symplectomorphisms.

\item {\rm [Monotonicity].} Let $Z'\subset Z\subset V\subset V'$. Then
$\CHZ(V,Z)\leq \CHZ(V',Z')$. In particular, $\CHZ(V,Z)\leq \CHZ(V)$.

\item {\rm[Homogeneity].} $\CHZ(V,Z,a\omega)=a \CHZ(V,Z,\omega)$, 
for any constant $a>0$.

\item {\rm [Normalization].} Assume that $M$ is a closed symplectic 
submanifold of a geometrically bounded symplectically aspherical manifold 
$W$ and let 
$U_r$ be a symplectic tubular neighborhood of $M$ 
in $W$ of radius $r>0$ (see Section \ref{sec:tubular} for the definition).
Then $\CHZ(U_r,M)=\pi r^2$.
\end{enumerate}
\end{Theorem}

Here all assertions, but the last one, are obvious.
Regarding the normalization assertion, we note that Theorem \ref{thm:main}
and Remark \ref{rmk:$C$} give the upper bound $\CHZ(U_r,M)\leq\pi r^2$. On the
other hand, it is straighforward to construct, for any positive 
$m< \pi r^2$,  a function $H\in \CH(U_r,M)$ with $\max H=m$ having no fast 
periodic orbits. This shows that $\CHZ(U_r,M)\geq \pi r^2$ and, thus, proves
the last assertion. 

A number of results on the Weinstein conjecture can be interpreted
as calculations of the relative Hofer--Zehnder capacity. 

\begin{Example}
\labell{exam:CHZ-c0}
Let $M$ be a compact symplectic manifold and 
$W=M\times\C^n$. Let $U_r=M\times B_r^{2n}$,
where $B_r^{2n}$ is the ball of radius $r>0$ in $\C^n$.
Then, as has been established in \cite{hv}, 
$$
\CHZ(U_r,M)=\pi r^2,
$$
provided that $r>0$ is small enough. If $\omega|_{\pi_2(M)}=0$, this holds
for any $r>0$. Note also that when $W$ is symplectically aspherical, this
can be proved similarly to the last assertion of Theorem \ref{thm:CHZ}.
Furthermore, when $\omega|_{\pi_2(M)}=0$, we also have
$\CHZ(U_r)=\pi r^2$ as is proved in \cite{fhv,ma}. Hence, in this case,
$\CHZ(U_r,M)=\pi r^2=\CHZ(U_r)$. We also refer the reader to
\cite{LT,lu} for further results in this direction.
\end{Example}

\begin{Example}
In general, the relative capacity $\CHZ(V,Z)$ need not be equal to the
capacity $\CHZ(V)$. For example, $\CHZ(V,Z) < \CHZ(V)$ whenever 
$\CHZ(V)<\infty$ and $Z$ is the closure of an open subset of $V$;
see Section \ref{sec:alm-ext} and, in particular, Proposition \ref{prop:ineq}.
Furthermore, it appears that the two capacities may differ for
small neighborhoods of a Lagrangian submanifold $Z$.
\end{Example}

\subsubsection{Almost existence}
As for the ordinary Hofer--Zehnder capacity, the finiteness of
relative capacity implies almost existence of periodic orbits:

\begin{Theorem}[Relative almost existence theorem]
\labell{thm:ae}
Assume that $\CHZ(V,Z)<\infty$ and let $H\colon V\to \R$ be a proper
smooth function constant on $Z$ and such that $H|_Z=\min H$. Then for
almost all (in the sense of measure theory) regular values $c$ in the range 
of $H$, the level $H=c$ carries a periodic orbit.
\end{Theorem}

The proof of this theorem, omitted here, is identical to the proof of the 
standard almost existence theorem, see  \cite[Section 4.2]{hz:book}. We
will further discuss the almost existence theorem and some elements of its
proof in Section \ref{sec:alm-ext}.

The notion of capacity considered above does not allow one to control the
homotopy class of periodic orbits. There are a number of ways to deal with
this shortcoming; see \cite{bps,sc}. Here we take an alternative approach
by restricting our attention to a collection of subsets of a fixed
symplectic manifold $V$. The resulting notion which henceforth we refer to
as the restricted relative Hofer--Zehnder capacity, even though lacking the 
flexibility of the capacities introduced in \cite{bps,sc}, is very simple 
and sufficient for our purpose to detect periodic orbits contractible in 
the ambient manifold $V$.

\subsubsection{Restricted Hofer--Zehnder capacity}
Let, as above, $(V,\omega)$ be a symplectic manifold,
$Z$ a compact subset of $V$, and $U$ an open
subset of $V$ containing $Z$.

\begin{Definition}
\labell{def:CHZW}
The \emph{restricted relative Hofer--Zehnder capacity}
$$
\CW(U,Z)\in (0,\infty]
$$ 
is defined as
$\sup\max H$, where the supremum is taken over $H\in \CH(U,Z)$ such that
all non-trivial, \emph{contractible in} $V$, periodic orbits of  $H$ 
are slow. When the dependence of the restricted capacity on $\omega$ is 
essential, we will use the notation $\CW(U,Z,\omega)$.
\end{Definition}

By definition, the restricted capacity is an invariant of $(U,Z)$ with 
respect to symplectomorphisms of the ambient manifold $V$. Furthermore,
it is clear that 
$$
\CHZ(U,Z)\leq\CW(U,Z)
$$
and
$$
\CHZ(U,Z)=\CW(U,Z)
\quad\text{if $U$ is simply connected.}
$$

Assume now that the ambient manifold is geometrically bounded and
symplectically aspherical. In accordance with our conventions, we denote it
by $W$.

\begin{Theorem}
\labell{thm:width}
Theorems \ref{thm:HZ-modify}, \ref{thm:CHZ}, and \ref{thm:ae} hold,
with obvious modifications (e.g., $V$ replaced by $U$), for the 
restricted relative Hofer--Zehnder capacity.
In particular, 
\begin{equation}
\labell{eq:added}
\CHZ(U_r,M)=\pi r^2= \CW(U_r,M), 
\end{equation}
when $M$ is a closed
symplectic submanifold of a geometrically bounded symplectically aspherical
manifold $W$. In the almost existence theorem, the flow of $H$ has
contractible in $W$ non-trivial periodic orbits on the level
$ H=c $ for almost all regular values $c$ in the range of $H$.
\end{Theorem}

This result can be easily verified by scrutinizing the proofs
of Theorems \ref{thm:HZ-modify}, \ref{thm:CHZ}, and \ref{thm:ae}. 
The normalization assertion, i.e., \eqref{eq:added}, readily follows from
the fact that Theorem \ref{thm:main} guarantees the existence of contractible 
periodic orbits.

As a consequence,  we have

\begin{Corollary}
\labell{cor:CHZW2}
Let $M$ be a closed symplectic 
submanifold of $W$ and let $H\colon W\to \R$ attain its absolute minimum 
$H=0$ on $M$. Then the levels $H=\eps$ carry contractible in $W$ periodic 
orbits for almost all small $\eps>0$. 
\end{Corollary}

This corollary generalizes a theorem of \cite{cgk} where the existence of 
periodic orbits is established for a dense set of small $\eps>0$. Note
also that another version of an almost existence theorem for convex manifolds
(also using the notion of relative capacity) has been recently proved 
in \cite{FS}.

\begin{Remark}
\labell{rmk:period}
As in \cite{hz:book}, it is easy to see that the constant
in Theorem \ref{thm:min} is in fact equal to $\CW(U,M)$ and there exists
a contractible orbit of period no greater than $\min_M H/\CW(U,M)$.
(When the orbits are not required to be contractible, one should replace
$\CW(U,M)$ by $\CHZ(U,M)$.)
\end{Remark}

\subsection{Application: the motion of a charge in a magnetic field}
\labell{subsec:magnetic}
Let $M$ be a closed Riemannian manifold and let $\sigma$ be a closed two-form
(magnetic field) on $M$. Equip $W=T^*M$ with the twisted symplectic 
structure $\omega=\omega_0+\pi^*\sigma$, where $\omega_0$ is the standard 
symplectic form on $T^*M$ and $\pi\colon T^*M\to M$ is the natural projection.
It is known that $(W,\omega)$ is geometrically bounded for any $\sigma$.
(We refer the reader to \cite{al,cgk,lu} 
for a discussion of this question.)
Finally, let $H$ be the standard kinetic energy Hamiltonian on $T^*M$.
The Hamiltonian flow of $H$ on $W$, called a \emph{twisted geodesic flow}, is
of interest because it describes, for example, the motion of a charge
on $M$ in the magnetic field $\sigma$.

In this setting, as a particular case of Corollary \ref{cor:CHZW2}, we have

\begin{Corollary}
\labell{cor:magn}
Assume that $(M, \sigma)$ is symplectically aspherical. Then for almost 
all small $\eps>0$, the energy level $ H=\eps $ carries a periodic orbit 
whose projection to $M$ is contractible.
\end{Corollary}

\begin{Example}
\labell{exam:horocycle}
Let $M$ be a closed surface with a metric of constant curvature equal
to $-1$ and let $\sigma$ be the area form on $M$. 
All orbits of the twisted geodesic flow on the levels $H=\eps$ are 
periodic and contractible in
$T^*M$ for $0<\eps<1/2$; for $\eps=1/2$, the level carries no periodic orbits;
and for $\eps>1/2$ the flow is smoothly equivalent to the geodesic flow
on the unit cotangent bundle $ST^*M$. 
In particular, the levels $H=\eps$ with $\eps>1/2$ carry
no periodic orbit with contractible projections to $M$. (We refer the
reader to, e.g., \cite{gi:survey} for proofs of these standard facts and 
further references.) As a consequence, we see that Corollary \ref{cor:magn}
does not extend to large values of $\eps$ as long as the orbits are
required to have contractible projections. Likewise, Theorem \ref{thm:main}
for contractible orbits
does not extend to arbitrarily large bounded neighborhoods of $M$.
Furthermore, set 
$W_\eps=\{H<\eps\}$. Now, it is easy to see from Theorem \ref{thm:main} that 
$\CHZ(W_\eps,M)=\CW(W_\eps,M)<\infty$ as long as $0<\eps<1/2$.
For $\eps>1/2$, we have $\CW(W_\eps,M)=\infty$ and it is not known
if $\CHZ(W_\eps,M)$ is still finite. 
\end{Example}

Corollary \ref{cor:magn} strengthens or complements a number of other results 
on periodic orbits of twisted geodesic flows. Under the hypotheses of the
corollary, the existence of  a dense set
of low energy levels with contractible periodic orbits has been proved
in \cite{cgk}. It is also known that when $\sigma$ is symplectic, but under
no assumptions on $\pi_2(M)$, there exists a sequence of energy values
$\eps_k\to 0$ such that the levels $ H=\eps_k$ carry contractible 
periodic orbits, \cite{gk2}. The same is true for any $\sigma\neq 0$, 
provided that $\pi_2(M)=0$, \cite{mac,pol2}. When $M$ is a torus, periodic
orbits (not necessarily contractible, e.g., if $\sigma=0$) exist on
almost all energy levels for any $\sigma$. In fact, in this case the
ordinary Hofer--Zehnder capacity of any bounded domain is finite; see,
e.g., \cite{gk1} and references therein, cf. \cite{ji}. This is also
true for the restricted Hofer--Zehnder capacity of small neighborhoods 
of $M$, when $\sigma\neq 0$ is exact and $\pi_2(M)=0$, \cite{FS}.

Furthermore, under suitable additional conditions, every low energy
level admits a number of periodic orbits when $\sigma$ is symplectic. 
We refer the reader to
\cite{gk1,ke} for some recent results and to, e.g., the survey
\cite{gi:survey} and references therein for a discussion of the results
obtained prior to 1995.

Corollary \ref{cor:magn} and the results quoted above suggest that
conjecturally, for any closed $\sigma$, every low energy level 
$ H=\eps$ carries a periodic orbit. However, this conjecture cannot
be established by purely symplectic topology methods: one should make use, 
in an essential way, of the fact that $H$ is convex along the fibers, as the 
the results of the next subsection indicate.

\subsection{Counterexamples}
\labell{sec:counterexamples}
We start with an example showing that in Corollary \ref{cor:CHZW2},
already for $W=\R^{2n}$, one cannot expect to have periodic orbits
on all levels of $H$ unless the minimum of $H$ is non-degenerate.
Note that in the non-degenerate case every low energy level $H=\eps$
carries at least $n$ distinct periodic orbits. This is the 
Weinstein--Moser theorem; see \cite{Mo:orbits,We:orbits}.

\begin{Proposition}
\labell{prop:non-exist-R2n}
For $2n\geq 6$, there exists a proper $C^\infty$-function 
$H\colon \R^{2n}\to\R$ whose only critical point is the origin,
where $H$ has an absolute minimum (say, $H(0)=0$), and
such that the levels $H=\eps_k$ carry no
periodic orbits for some sequence of regular values $\eps_k\to 0$. 
The levels $H=\eps_k$ are isotopic to the standard sphere in $\R^{2n}$.
When $2n=4$, there is a $C^2$-smooth function with these properties.
\end{Proposition}

The construction of $H$ as well as the constructions of the other
two counterexamples below will be given in Section 
\ref{sec:counterexamples-prfs}.

\begin{Remark}
The function $H$ has a very degenerate (flat) minimum at the origin.
One can also modify the construction of $H$ so that $H$ has two
critical sets: a non-degenerate minimum $H(0)=0$ and a critical set
$H=\eps>0$, diffeomorphic to $S^{2n-1}$, and $\eps_k\to\eps$.
 These results were stated in \cite{gi:barcelona} without proof. 
\end{Remark}

The next counterexample concerns Hamiltonian flows on $T^*S^n$ with the 
standard symplectic structure.

\begin{Proposition}
\labell{prop:non-exist}
For $n\geq 3$, there exists a $C^\infty$-function
$H\colon T^*S^n\to \R$ attaining its absolute minimum $H=0$ at the zero section
and such that for some sequence of regular values $\eps_k\to 0$ the levels 
$ H=\eps_k$ carry no periodic orbits of the Hamiltonian flow of $H$
with respect to the standard symplectic structure $\omega_0$ on
$T^*S^n$. The zero section is the only critical set of $H$ and the
levels $H=\eps_k$ are isotopic to the unit cotangent bundle $ST^*S^n$.
When $n=2$, there exists a $C^2$-smooth function with these
properties. Moreover, the same is true for some exact non-zero magnetic
field $\sigma$ on $S^n$, $n\geq 2$.
\end{Proposition}

This proposition shows that in the results of \cite{hv-ct,vi} concerning
the Weinstein conjecture in the cotangent bundles the contact type
condition cannot, in general, be omitted.

Finally, let us turn to magnetic flows for non-degenerate magnetic fields.
Let $\sigma$ be the standard symplectic form on $M=\CP^n$ and, as above,
$\omega=\omega_0+\pi^*\sigma$ the twisted symplectic form on $W=T^*\CP^n$.

\begin{Proposition}
\labell{prop:non-exist-CPn}
The assertions of Proposition \ref{prop:non-exist} hold when
$(T^*S^n,\omega_0)$ is replaced by $(T^*\CP^n,\omega)$ with a 
$C^\infty$-smooth Hamiltonian $H$ for $n\geq 2$ and a $C^2$-smooth Hamiltonian 
for $n=1$.
\end{Proposition}

\begin{Remark}
In this example, $M$ is not symplectically aspherical. It appears plausible 
that similar examples exist for symplectically aspherical manifolds.
\end{Remark}

\begin{Remark}
As in Proposition \ref{prop:non-exist-R2n}, the function $H$ in 
Propositions \ref{prop:non-exist} and \ref{prop:non-exist-CPn}
has a very degenerate (flat) minimum along the zero section.
\end{Remark}

\section{Time-dependent Hamiltonian flows}
\labell{sec:time-dep}

So far we have been concentrating on the existence problem for
non-trivial periodic orbits of autonomous Hamiltonians. The main technical
tool in our dealing with this problem has been Theorem \ref{thm:main} which
guarantees that an independent of time Hamiltonian $H\in \CH(W,M)$ 
with sufficiently large variation has a non-trivial 
contractible one-periodic orbit.

An analogue of this result holds for periodic in time Hamiltonians. In this
setting, inspired by Arnold's conjecture, the requirement that the orbit
be non-trivial appears to make little sense. Furthermore, once this
requirement is dropped, the existence result can be established under
conditions much less restrictive than in Theorem \ref{thm:main}. Below
we prove only the simplest version of the theorem one can expect to 
hold for time-dependent Hamiltonians (see Remark \ref{rmk:$H>0$}), 
for our goal is to emphasize
the difference between the autonomous and time-dependent cases rather than
to analyze the time-dependent case in detail.

Let, as in Section \ref{sec:main}, $W$ be a geometrically bounded 
symplectically aspherical manifold and let $M$ be a closed symplectic 
submanifold 
of $W$. Consider a one-periodic in time Hamiltonian $H\colon S^1\times W\to\R$
supported in $S^1\times U$, where $U$ is an open set containing $M$.
Then, along the lines of \cite{bps}, we have the following analogue of 
Theorem \ref{thm:main}:

\begin{Theorem}
\labell{thm:time-dependent}
Assume that $U$ is a sufficiently small neighborhood of $M$. Then 
for every $H\geq 0$ supported
in $S^1\times U$ and such that $\min_{S^1\times M} H> 0$, the Hamiltonian 
flow of $H$ has a contractible one-periodic orbit with positive action.
If the periodic orbits of $H$ with positive action are non-degenerate, the 
number of these orbits is no less than the sum of Betti numbers of $M$ 
with $\Z_2$ coefficients.
\end{Theorem}

This theorem is an immediate consequence of the following

\begin{Theorem}
\labell{thm:time-dependent2}
Assume that $U$ is a sufficiently small neighborhood of $M$ and let
$H\geq 0$, supported in $S^1\times U$, be such that $\min_{S^1\times M} H> 0$.
Then, for every small $a>0$, there exists an epimorphism 
$\HF^{[a,\,\infty)}(H)\to H_*(M;\Z_2)$.
\end{Theorem}

Theorem \ref{thm:time-dependent2} will be proved in Section \ref{sec:pfs}.

Theorem \ref{thm:time-dependent} holds under much weaker hypotheses
than its autonomous counterpart, Theorem \ref{thm:main}. However, the trade-off
is that Theorem \ref{thm:time-dependent} cannot guarantee the existence
of non-trivial periodic orbits. As is well known, 
for any $C^2$-small autonomous Hamiltonian $H$,
one-periodic orbits of the flow of $H$ are trivial 
(cf. the construction in Example \ref{exam:actions} with $c>0$ small).
Moreover, in sharp contrast with Theorem \ref{thm:main}, 
even when $M$ is a point in the setting of Theorem \ref{thm:time-dependent}, 
there exists a Hamiltonian $H$ with arbitrarily large 
$\min_{S^1\times M} H$ such that all its positive actions are
arbitrarily small (see Example \ref{exam:polt}).\footnote{The authors 
are greateful to Leonid Polterovich for this remark.} In other words,  
the value of $a$ in Theorem \ref{thm:time-dependent2} cannot be
expressed in terms of $\min_{S^1\times M} H$ only. Furthermore, Example
\ref{exam:polt} also shows that there exists a 
compactly supported Hamiltonian $H$ (not satisfying
$H\geq 0$) such that $\min_{S^1\times M} H$ is arbitrarily large
but all the action values for $H$ are non-positive and the time-one
flow of $H$ is non-trivial. These examples demonstrate
that for a small neighborhood of $M$ the relative symplectic capacities 
introduced in \cite{bps} are trivial.

Further comparing Theorems \ref{thm:time-dependent} and 
\ref{thm:time-dependent2} with the results of \cite{bps} note that
when $M$ is a Lagrangian submanifold and the orbits are
sought in a non-trivial homotopy class $\alpha$, the assumption
$\min_{S^1\times M}H> C_\alpha>0$ from \cite{bps} is clearly logical
and necessary. However, when the orbits looked for are contractible
the role of the condition $\min_{S^1\times M}H>0$ is less apparent, 
in particular, if the requirement that the action is positive replaced by 
that the action is non-zero. 

More specifically, it is possible that the first assertion of Theorem
\ref{thm:time-dependent} still holds when $H\geq 0$, but the condition
that $\min_{S^1\times M} H>0$ is replaced by the requirement that $H$
is not identically zero. When $W$ is convex, this can be proved using
the monotonicity of the Schwarz action selector
(cf. \cite{FS,sc}).  When $H$ is not assumed to be non-negative, the
assumption that the time-one flow of $H$ is non-trivial should be
added and in this case one should look for the orbits with non-zero
action. This would lead to a partial generalization of theorem of
C. Viterbo, \cite{Vi:thm}.  (Note that when $H\geq 0$ is not
identically zero, the time-one flow is non-trivial, at least when the
symplectic manifold is exact, since the Calabi invariant of the flow
is positive; see \cite{McDS}.) For convex manifolds, this version of
Viterbo's theorem has been proved in \cite{FS}. However, the methods
of \cite{FS} heavily rely on a somewhat different definition of Floer
homology, making use of convexity and resulting in the homology
defined for all intervals of actions. This approach does not readily
extend to geometrically bounded symplectic manifolds.

\section{Preliminaries}
\labell{sec:prelim}

\subsection{Symplectic tubular neighborhoods}
\labell{sec:tubular}
Let $M$ be a compact symplectic submanifold of a symplectic manifold 
$(V,\omega)$.
Denote by $E\to M$ the symplectic normal bundle to $M$. Recall that a 
neighborhood of the zero section in $E$ has a natural symplectic structure 
$\omega_E$. Moreover, 
on this neighborhood, there exists a fiberwise quadratic Hamiltonian
whose flow is periodic. This can be seen as follows.

Let us equip $E$ with a
Hermitian metric compatible with the fiberwise symplectic structure
on $E$ and denote by $\rho\colon E\to \R$ the square of the fiberwise 
norm, i.e., $\rho(z)=\|z\|^2$. Recall that $E$ has a canonical 
fiberwise one-form whose
differential is the fiberwise symplectic form. (The value of this form
at $z\in E$ is equal to the contraction of the fiberwise symplectic form with 
$z$.) Fixing a Hermitian connection on $E$, we extend this fiberwise one-form
to a genuine one-form $\theta$ on $E$. 

Then the form
$$
\omega_E=\12 d\theta+\sigma
$$
is symplectic on a neighborhood of the zero section in $E$. Here we 
have identified $\sigma=\omega|_M$ with its pull-back to $E$.

By the symplectic neighborhood theorem, the neighborhood 
$\{\rho< r^2\}$ of $M$ in $(E,\omega_E)$
is symplectomorphic to a neighborhood of $M$ in $(V,\omega)$
for some $r>0$. From now on, we denote this neighborhood by $U_r$,
assume the identification of $\{\rho<r^2\}$ and $U_r$,
and refer to $U_r$ together with this identification
as a \emph{symplectic tubular neighborhood of $M$ in $V$ of radius
$r$}. In particular, in what follows, $\omega=\omega_E$  
and $\rho$ is regarded as a function on a neighborhood of $M$ in $V$. 
Sometimes, we will write $\rho(x)$ as $\| x \|^2$.

It is not hard to see that all orbits of the Hamiltonian flow of the 
function $\12 \rho\colon E\to\R$
are periodic with period $2\pi$, just as for the square of the standard norm
on $\R^{2n}$. This fact will be essential for the calculation of the Floer
homology of a small tubular neighborhood of $M$ in $V$.

\subsection{Floer homology}
\labell{sec:floer}
In this section, we will recall a few facts concerning Floer homology needed
for the proofs.  The reader interested in a detailed 
treatment of this material should consult the Floer's papers 
\cite{fl1,fl2,fl3} or, for example,  \cite{hz:book,sa} for 
a general introduction to Floer homology, and \cite{bps,cfh,cfhw,cgk,fh,fhw}
for the definition and properties of symplectic homology.

\subsubsection{The definition of filtered Floer homology}
Let, as in Section \ref{sec:main}, $(W,\omega)$ be a geometrically 
bounded symplectically aspherical manifold. 
Denote by $\mathcal{H}$ the space of smooth, compactly 
supported Hamiltonians $H\colon S^1 \times W\to \R$. To each 
$H \in \mathcal{H}$, we associate the action functional $A_H$, defined
by \eqref{eq:action}, on the space of smooth contractible loops in $W$.
The critical points $\mathcal{P}(H)$ of $A_H$ are exactly contractible 
one-periodic orbits of the time-dependent Hamiltonian flow of $H$. 
The set of critical values of $A_H$ is called the \emph{action spectrum} 
of $H$ and we denote it by
$$
\mathcal{S}(H)=\{ A_H(x) \mid x \in \mathcal{P}(H) \}.
$$
It is known that $\mathcal{S}(H)$ is compact and nowhere dense in $\R$,
\cite{hz:book,sc}.

The Floer homology of $H$ for a certain interval of actions 
is the homology of the (relative) Morse complex of  $A_H$ on the space of all
contractible loops.
However, when $W$ is not compact, every point in the complement of
$\supp H$ is a degenerate one-periodic orbit, i.e., a critical point of
$A_H$, with zero action. To avoid this set, we will only
consider the homology generated by the contractible one-periodic
orbits with action in an interval that does not contain zero.

To make this description more precise, let us first recall, following
\cite{cfh}, 
the definition of the Floer homology $\HF^{[a,\,b)}(H)$ for the negative
range of actions $a<b<0$. For a fixed $a<0$ such that 
$a\not\in\mathcal{S}(H)$, let
$$
\mathcal{P}^a(H)=\{x \in \mathcal{P}(H) \mid A_H(x) <a\}.
$$
Assume first that $H$ satisfies the following condition:
\begin{equation}
\labell{eq:cond}
\text{Every one-periodic orbit $x \in \mathcal{P}^a(H)$ is nondegenerate.}
\end{equation}
Since $c_1(TM)|_{\pi_2(M)}=0$, the elements of $\mathcal{P}^a(H)$
are graded by the Conley-Zehnder index $\mu_{\CZ}$ (see, e.g., \cite{sa})
and the Floer complex of $H$ for actions less than $a$ is
the graded $\Z_2$-vector space
$$
\CF^a(H) = \bigoplus_{x \in {\Pp^a (H)}} \Z_2 x.
$$
Note that for the action functional \eqref{eq:action}, 
the Hamiltonian vector field $X_H$ is given by $dH = - i_{X_H} \omega$
(differing by sign from the Hamiltonian vector field in, say, \cite{sa})
and the direction of $X_H$ effects the sign of
Conley--Zehnder indices.

To define the Floer boundary operator, we first fix an almost
complex structure $J_{gb}$ for which $(W, \omega)$ is
geometrically bounded as in \cite{cgk}. Let $\Jj$ be
the set of smooth $t$-dependent $\omega$-tame almost complex
structures which are $\omega$-compatible near $\supp(H)$ and are equal
to $J_{gb}$ outside some compact set. Each $J \in \Jj$ gives rise to a
positive--definite bilinear form on the space of contractible loops in
$W$. We can then consider
the moduli space $\Mm (x,y,H,J)$ of \emph{downward} gradient-like
trajectories of $A_H$ which go from $x$ to $y$ and have finite
energy. For a dense subset, $\Jj_{reg}(H) \subset \Jj$, each moduli space
$\Mm (x,y,H,J)$ is a smooth manifold of dimension $\mu_{\CZ}(x)
-\mu_{\CZ}(y)$.

As usual, the Floer boundary operator is then defined by
\begin{equation}
\labell{eq:floer-d}
\partial^{H,J} x = \sum_{y \in \mathcal{P}^a(H) \text{ with } 
\mu_{\CZ}(x) -\mu_{\CZ}(y) =1} \tau(x,y)y,
\end{equation}
where $\tau(x,y)$ stands for the number (mod $2$) of elements in $\Mm
(x,y,H,J)/ \R$ and $\R$ acts (freely) by translation on the
gradient-like trajectories. The operator $\partial^{H,J}$ 
satisfies $\partial^{H,J} \circ
\partial^{H,J} =0$ and the resulting Floer homology groups
$\HF^a(H)$ are independent of the choice of $J \in \Jj_{reg} (H)$.

\begin{Remark}
\labell{rem:compact1} Since $(W, \omega)$ with $J_{gb}$ is
geometrically bounded and $H$ is compactly supported, there is a
uniform $C^0$-bound for the elements of $\Mm (x,y,H,J)$ (see, for
example, Chapter V in \cite{al}). Hence, the compactness of the
appropriate moduli spaces follows by the usual arguments.
It is unclear whether or not $\HF^a(H)$ depends on the choice of
$J_{gb}$. It is independent of this choice if the set of almost
complex structures for which $W$ is geometrically bounded is
connected.
\end{Remark}

For any pair $a<b$, set
$$
\mathcal{H} ^{a,b} = \{ H\in \mathcal{H} \mid a,b \notin \mathcal{S} (H) \}.
$$
Assume  that $H$ has property \eqref{eq:cond} for $b$ (and hence for $a$). 
(Note that this can only happen when $a<b<0$.) Then the complexes
$\CF^a(H)$ and $\CF^b(H)$ are defined and 
$\CF^a(H)$ is a subcomplex of $\CF^b(H)$. By definition, $\HF^{[a,\,b)}(H)$ is
the homology of the quotient complex $\CF^{[a,\,b)}(H) = \CF^b
(H)/ \CF^a (H)$ with the induced boundary operator.

The set $\mathcal{H}^{a,b}$ is open in $\mathcal{H}$ with
respect to the strong Whitney $C^{\infty}$-topology. Moreover, in
each component of $\mathcal{H}^{a,b}$, the functions with the
property \eqref{eq:cond} holding for $b$ form a dense set.
For any function $H\in \mathcal{H}^{a,b}$, we define the Floer
homology $\HF^{[a,\,b)}(H)$ as $\HF^{[a,\,b)}(K)$ where $K$ is a small
perturbation of $H$ such that $K\in \mathcal{H}^{a,b}$ and $K$ satisfies
\eqref{eq:cond} for $b$. A version of Floer's continuation then shows that
$\HF^{[a,\,b)}(K)$ is independent of $K$ as long as $K$ is close to $H$. Hence,
$\HF^{[a,\,b)}(H)$ is well defined.

When $0<a<b$, the condition \eqref{eq:cond} for $a$ or $b$ is never
satisfied and we adopt
a different, and more naive, approach to the definition of Floer homology.
Namely, we simply work with the complex generated by periodic orbits with
action in $(a,\,b)$. To be more precise, assume that 
$H\in \mathcal{H} ^{a,b}$ and let
$$
\mathcal{P}^{a,b}(H)=\{x \in \mathcal{P}(H) \mid a<A_H(x) <b\}.
$$
Also, let us temporarily assume that the following condition holds:
\begin{equation}
\labell{eq:cond2}
\text{Every one-periodic orbit $x \in \mathcal{P}^{a,b}(H)$ 
is nondegenerate.}
\end{equation}
Consider the $\Z_2$-vector space 
$$
\CF^{[a,\,b)}(H) = \bigoplus_{x \in {\Pp^{a,b} (H)}} \Z_2 x
$$
graded by the Conley--Zehnder index. The differential $\partial^{H,J}$
is defined similarly to \eqref{eq:floer-d}, but with summation extending
only to $y\in \Pp^{a,b}$. As before, $\left(\partial^{H,J}\right)^2=0$,
and we set $\HF^{[a,\,b)}(H)$ to be the cohomology of the resulting complex.
Since the functions satisfying \eqref{eq:cond2} for $a<b$ are dense in
$\mathcal{H}^{a,b}$, a small perturbation argument as above allows us to
define  $\HF^{[a,\,b)}(H)$ for any $H\in \mathcal{H} ^{a,b}$.

Note that this construction can also be used when $a<b<0$ and in this case
the two definitions lead to the same complex $\CF^{[a,\,b)}(H)$.

The behavior of Floer homology when the interval of actions is shrunk is
described by a long exact sequence. Namely, assume that $a<b<c$,
none of these points is in $\mathcal{S}(H)$, and $0\not\in [a,\,c]$. Then
we have the exact sequence 
$$
\ldots\to\HF^{[a,\,b)}_{*}(H)\to
\HF^{[a,\,c)}_{*}(H)\to
\HF^{[b,\,c)}_{*}(H)\to
\HF^{[a,\,b)}_{*-1}(H)
\to\ldots .
$$
Indeed, for either of the above definitions of the Floer complex,
we obviously have the exact sequence of complexes
$$
0\to\CF^{[a,\,b)}(H)\to\CF^{[a,\,c)}(H)\to\CF^{[b,\,c)}(H)\to 0,
$$
which induces the required exact sequence in Floer homology.

In particular, this shows that the end points of the interval $[a,\,b]$
can be continuously varied without changing $\HF^{[a,\,b)}(H)$ as long
as $a$ and $b$ stay away from $\mathcal{S}(H)$.

\begin{Remark}
We will use the Floer homology $\HF^{[a,\,b)}(H)$ with $0<a<b$ to prove that
the function $H$ from Theorem \ref{thm:main} has a non-trivial one-periodic
orbit with positive action. In this setting, one can easily avoid making use
of the above construction by replacing $H$ by $-H$ or by working with
cohomology, i.e., considering the differential defined by counting upward
gradient-like trajectories. However, we have found the setting of positive 
functions and downward trajectories visually more pleasing, which has 
motivated our choice.
\end{Remark}

\subsubsection{Monotone homotopies and monotone homotopy 
invariance of Floer homology}
Let $H,K \in \mathcal{H}^{a,b}$ be two functions with $H(t,x) \geq
K(t,x)$ for all $(t,x) \in S^1 \times W$. Then there exists a monotone
homotopy $s \mapsto K_s$ from $H$ to $K$, i.e., a family of functions $K_s$
such that $s\mapsto K_s(t,x) $ is monotone decreasing for all $(t,x)$ and
$$
K_s=\begin{cases}
H &\text{for $s\in (-\infty,-1]$,}\\
K & \text{for $s\in [ 1, \infty)$.}
\end{cases}
$$
Such a homotopy induces a Floer chain map
$$
\CF^{[a,\,b)}(H)\to \CF^{[a,\,b)}(K),
$$ 
and hence a homomorphism of Floer homology
$$
\sigma_{KH} \colon \HF^{[a,\,b)}(H) \ra \HF^{[a,\,b)}(K).
$$
Note that $K_s$ is not required to be in $\mathcal{H}^{a,b}$.

The following facts concerning these homomorphisms are well known;
see, e.g., \cite{cfh,fh,vi} and \cite[Sections 4.4 and 4.5]{bps}:
\emph{The homomorphism $\sigma_{KH}$ is independent of the
choice of the monotone homotopy $K_s$ and has the following
properties:}
\begin{align*}
\sigma_{KH} \circ \sigma_{HG}& = \sigma_{KG} \text{ for } G \geq H \geq K,\\
\sigma_{HH} &= \id  \text{ for every $H \in  \mathcal{H}^{a,b}.$}
\end{align*}
Furthermore, Floer homology is homotopy invariant in the following sense:
\emph{Assume that $K_s \in \mathcal{H}^{a,b}$ for all $s \in [0,1]$.
Then $\sigma_{KH}$ is an isomorphism.}

This shows that the only way in which the map
$\sigma_{KH}$ can fail to be an isomorphism is if periodic orbits,
with action equal to $a$ or $b$, are created during the homotopy.

\subsubsection{Calculations of Floer homology}
The main tool used in this paper to calculate Floer homology
is a theorem of Po\'zniak, \cite{poz}, which equates filtered 
Floer homology and the ordinary homology of a connected Morse--Bott 
non-degenerate set of periodic orbits.

Recall that a subset $P \subset \mathcal{P}(H)$ is said to be a
\emph{Morse--Bott non-degenerate manifold of periodic orbits} 
if the set $C_0 =
\{x(0) \mid x \in P\}$ is a compact submanifold of $W$ and
$T_{x_0}C_0 = \ker (D \phi_H^1(x_0) - \id)$ for every $x_0 \in
C_0.$ Here $\phi_H^1$ is the time-one flow of $X_H$.

For such sets of periodic orbits we have the following result
which holds for geometrically bounded, symplectically aspherical 
manifolds.

\begin{Theorem}
\labell{thm:mb} {\rm (Po\'zniak, \cite[Corollary 3.5.4]{poz})}
Let $ a < b $ be outside of the action spectrum of $H$
and such that $[a,\, b]$ does not contain zero. Also, suppose
that the set $P = \{ x \in \Pp (H) \, \mid \, a < \Aa_H < b\}$ is
a connected Morse--Bott manifold of periodic orbits. Then
$\HF^{[a,\,b)}(H)$ is isomorphic to  $H_* (P; \Z_2)$.
\end{Theorem}

Note that this isomorphism does not preserve the grading, i.e.,
$\HF^{[a,\,b)}_{*}(H)=H_{*-s} (P; \Z_2)$, where the shift $s$ depends
on the behavior of $H$ near $P$. 

We will also need the following elementary (essentially, trivial) 
observation which sometimes allows one to
extend Po\'zniak's isomorphism in a particular degree to the case where
$P$ is disconnected.

\begin{Lemma}
\labell{lemma:poz}
Let $ a < \gamma < b $ be outside of the action spectrum of $H$
and such that $[a,\, b]$ does not contain zero.
\begin{enumerate}
\item Suppose that 
$$
\HF^{[\gamma,\,b)}_{n_0+1}(H)=\HF^{[\gamma,\,b)}_{n_0}(H)=0.
$$
Then the natural map $\HF^{[a,\,\gamma)}_{n_0}(H)\to\HF^{[a,\,b)}_{n_0}(H)$
is an isomorphism.

\item Suppose that 
$$
\HF^{[a,\,\gamma)}_{n_0-1}(H)=\HF^{[a,\,\gamma)}_{n_0}(H)=0.
$$
Then the natural map $\HF^{[a,\,b)}_{n_0}(H)\to\HF^{[\gamma,\,b)}_{n_0}(H)$
is an isomorphism.
\end{enumerate}
\end{Lemma}

We will use this lemma in the situation where Po\'zniak's theorem
applies to the intervals $[\gamma,\,b]$ or $[a,\,\gamma]$. Then the Floer 
homology groups for these intervals vanish in degrees outside of a certain 
range $[\mu_-,\,\mu_+]$ of Conley--Zehnder indices. Therefore,
$\HF^{[a,\,\gamma)}_{n_0}(H)\to\HF^{[a,\,b)}_{n_0}(H)$ is an isomorphism
when, for instance, $n_0> \mu_+$ and 
$\HF^{[a,\,b)}_{n_0}(H)\to\HF^{[\gamma,\,b)}_{n_0}(H)$ is an isomorphism
when $n_0-1 > \mu_+$. This line of reasoning has been used in \cite{cgk}
to calculate the Floer homology in a given degree for large intervals
of actions.

\begin{proof}[Proof of Lemma \ref{lemma:poz}]
To prove the first assertion, consider the exact sequence
$$
\HF^{[\gamma,\,b)}_{n_0+1}(H)\to
\HF^{[a,\,\gamma)}_{n_0}(H)\to
\HF^{[a,\,b)}_{n_0}(H)\to
\HF^{[\gamma,\,b)}_{n_0}(K_+).
$$
The first and the last group in this sequence vanish. Hence, the middle 
map in the exact sequence is an isomorphism.

In a similar vein, to prove the second assertion, we consider the exact 
sequence
$$
\HF^{[a,\,\gamma)}_{n_0}(H)\to
\HF^{[a,\,b)}_{n_0}(H)\to
\HF^{[\gamma,\,b)}_{n_0}(H)\to
\HF^{[a,\,\gamma)}_{n_0-1}(H).
$$
Here again the first and the last group vanish, and therefore
the middle map is an isomorphism.
\end{proof}

\section{Proofs of Theorems \ref{thm:main} and \ref{thm:time-dependent2}} 
\labell{sec:pfs}

\subsection{Proof of Theorem \ref{thm:main}}
First note that it suffices to prove the theorem 
in the case where $U$ is a tubular symplectic neighborhood, say $U_R$, of 
$M$ in $W$. Then the theorem and Remark \ref{rmk:$C$} are consequences of the 
following result giving the exact value of the constant $C$ for $U_R$.

\begin{Theorem}
\labell{thm:cap-nbd}
Assume that $M$ is a closed symplectic submanifold of a geometrically bounded
symplectically aspherical manifold $W$. Then for $R>0$ small enough and 
for every $H\in \CH(U_R,M)$ with $\max H> \pi R^2$, the Hamiltonian flow 
of $H$ has a non-trivial contractible one-periodic orbit with action 
in the interval $(\max H, \max H+ \pi R^2]$.
\end{Theorem}

This theorem, in turn, is based on

\begin{Proposition} 
\labell{prop:floer}
Let $H\in \CH(U_R,M)$ and $\max H > \pi R^2$. Then
there exist functions $K_+$ and $K_-$, supported in $U_R$, such that
$K_-<H<K_+$ and 
$$
\Z_2=\HF^{[a,\,b)}_{n_0}(K_+)\to \HF^{[a,\,b)}_{n_0}(K_-)=\Z_2
$$
is an isomorphism for some constants $\max H<a<b$, not in the action
spectra of $K_\pm$, and $n_0=\frac{1}{2}(\codim M-\dim M)+1$. The functions
$K_\pm$ can be chosen so that $b$ is arbitrarily close to $\max H+\pi R^2$.
\end{Proposition}

To prove Theorem \ref{thm:cap-nbd}, we factor the monotonicity isomorphism
as 
$$
\Z_2=\HF^{[a,\,b)}_{n_0}(K_+)\to \HF^{[a,\,b)}_{n_0}(H)\to
\HF^{[a,\,b)}_{n_0}(K_-)=\Z_2 .
$$
It follows that $\HF^{[a,\,b)}_{n_0}(H)\neq 0$. Thus, $H$
must have a contractible one-periodic orbit with action in the
interval $[a,\,b)$. Since $a>\max H$, this orbit is non-trivial. 

To complete the proof of Theorems \ref{thm:main} and \ref{thm:cap-nbd}, 
it remains to prove the proposition.

\subsection{Proof of Proposition \ref{prop:floer}}~
The idea of the proof is to pick $K_-$ and $K_+$ depending only on $\rho$
and squeezing $H$ from above and below as tightly as possible
(see Fig.\ 1). This will
guarantee that $[a,\, b)$ with the required properties does exist. This
interval contains more than one value of the action spectrum of $K_+$ or
$K_-$. However, only one one-periodic level of $K_\pm$ contributes 
$\Z_2$ in degree $n_0$ to the Floer homology of $K_\pm$ and the 
remaining levels with actions in $[a,\, b)$
contribute to the Floer homology in degrees either less than $n_0-1$ or 
greater than $n_0+1$. Then the exact sequence argument (Lemma \ref{lemma:poz})
shows that the interval $[a,\, b)$ can be shrunk to contain only the action 
value essential in degree $n_0$ without changing the homology in this
degree. Applying Po\'zniak's theorem, we see that 
$\HF_{n_0}^{[a,\,b)}(K_\pm)=\Z_2$. A similar argument shows that the Floer 
homology in degree $n_0$ remains constant in the course of a monotone
homotopy from $K_+$ to $K_-$ even though $a$ and $b$ do not stay outside
of the action spectrum.

\subsubsection{The definitions of $K_\pm$} The graphs of functions $K_\pm$ 
are shown in Fig.\ 1.\footnote{Here we break a recent but 
well-established tradition
to define the functions explicitly and in every detail, as is done for 
example in \cite{bps,cgk},
and revert to describing only the essential features of the functions, cf.
\cite{fhw}.}
 These functions depend only on $\rho$ and in what
follows we do not distinguish, for the sake of brevity,  the functions of 
$\rho$ from the corresponding functions on $U_R$. The shape of the functions
is similar to that used in \cite{bps}, however $\max K_\pm$ are chosen so
that these functions bound $H$ from below and above as tightly as possible.
Let us now specify some details in the definitions of $K_\pm$.

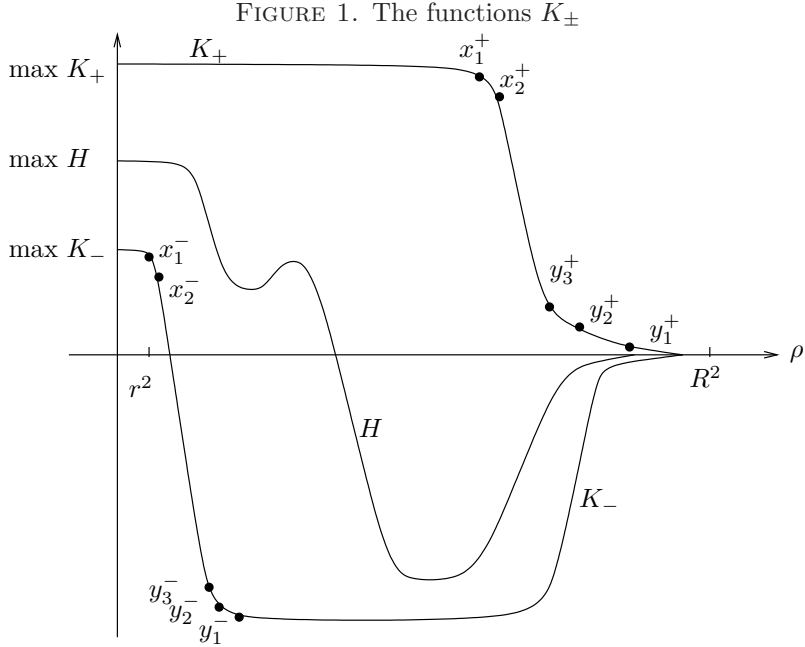
\begin{figure}[htbp]
\begin{center}
\caption{The functions $K_\pm$}
\input{functions.pstex_t}
\label{figure:functions}
\end{center}
\end{figure}

The function $K_+$ is constant and equal to $\max K_+$ until $\rho$ 
becomes nearly equal to $R^2$. Then the function rapidly decreases to zero
and is
identially zero when $\rho$ is very close to $R^2$. The slope of $K_+$, on the
interval where this function is non-constant linear, is not an integer 
multiple of $\pi$. Hence, the Hamiltonian flow of
$K_+$ has non-trivial one-periodic orbits on a finite sequence of levels 
where $\rho$ assumes values:
$$
x^+_1 < x^+_2 < x^+_3 < \ldots < \ldots y^+_3  < y^+_2 < y^+_1 < R^2.
$$
The points $x^+_l$ are located where the value of the function is
still close to $\max K_+$ and the points $y^+_l$ are located where
the value of the function is close to zero. 
Note that $x^+_1\approx R^2$ due to our choice of $K_+$.

The function $K_-$ is constant and equal to $\max K_-$ until $\rho$ 
becomes nearly equal to $r^2$ for a sufficiently small constant $r$
such that $0<r<R$, to be specified later. Then the function rapidly decreases 
to $\min K_-$. The value $\min K_-$ is chosen so that $K_-<H$. Hence,
$\min K_-$ is negative if $H$ assumes negative values and we can
take $\min K_-=0$ if $H\geq 0$. In what follows, we describe $K_-$ in
the former case. The function $K_-$ remains constant and equal to
$\min K_-$ until $\rho$ nearly reaches $R^2$. Then the function rapidly
increases to zero and becomes identically zero for $\rho$ very close to
$R^2$. On the intervals where $K_-$ is non-constant linear, the slopes are 
chosen not to be integer multiples of $\pi$. 

The Hamiltonian flow of $K_-$ has non-trivial one-periodic orbits on four
finite sequences of levels. The first two of them,
$$
r^2 < x^-_1 < x^-_2 < x^-_3 < \ldots < \ldots < y^-_3 < y^-_2 < y^-_1,
$$
are located where $K_-$ is decreasing and 
$K_-\approx \max K_-$ and $K_-\approx \min K_-$, respectively. 
Note that $y^-_1\approx r^2$ by the construction of $K_-$. The other
two sequences of levels are located where $K_-$ is increasing. For these
levels the actions are negative and hence the periodic orbits on these
levels do not contribute to the Floer homology $\HF^{[a,\,b)}(K_-)$.

Note that the points $x^\pm_l$ and $y^\pm_l$ are labelled so that
the slope of $K_\pm$ increases from $x^\pm_l$ to $x^\pm_{l+1}$
and from $y^\pm_l$ to $y^\pm_{l+1}$, and 
the periodic orbits on the levels $\rho=x^\pm_l$ and $\rho=y^\pm_l$ have
multiplicity $l$.

Particular attention should be given to the choice of $\max K_\pm$.
To describe how these maximal values are chosen,
denote by $A(x^\pm_l)$ and $A(y^\pm_l)$ the action of $K_\pm$ 
on the periodic orbits occuring on the levels $\rho=x^\pm_l$ and
$\rho=y^\pm_l$, respectively. Then we require that
\begin{equation}
\labell{eq:max}
\max K_-< \max H < \max K_+ < A(x^-_1) < A(x^+_1).
\end{equation}

Let us show that this choice is possible. The value $\max H$ is given
and we also know that $H$ is constant near $M$. Then we can chose $K_-$
and $r>0$ so that $K_-< H$ near $M$ and 
$$
\max K_-< \max H < A(x^-_1)\approx \max K_- + \pi r^2.
$$
Note that this can be done for an arbitrarily small $r>0$.
Finally, pick $\max K_+$ so that
$$
\max H < \max K_+ < A(x^-_1) < A(x^+_1)\approx \max K_+ + \pi R^2.
$$
This is clearly possible since $r>0$ is small and $R$ is fixed.

Finally note that the functions $K_\pm$ are strictly convex or
concave outside of the intervals where these functions are constant or
linear. This ensures that the energy levels $\rho=x^\pm_l$ and $\rho=y^\pm_l$
are Morse--Bott non-degenerate. Furthermore,
the functions $K_\pm$ can be chosen so that all action values
$A(x^\pm_l)$ and $A(y^\pm_l)$ are distinct.

\subsubsection{Periodic orbits, actions and Conley--Zehnder indices for 
$K_\pm$} When $[\alpha,\,\beta)$ is an interval of (positive) actions
containing only one of the points $A(x^\pm_l)$ and $A(y^\pm_l)$, the
Floer homology $\HF^{[\alpha,\,\beta)}(K_\pm)$ can be determined by 
Po\'zniak's theorem, \cite{poz}, (see Theorem \ref{thm:mb}).
Namely, we have
\begin{equation}
\labell{eq:poz}
\HF^{[\alpha,\,\beta)}_*(K_\pm)=H_{*-s}(SM;\Z_2),
\end{equation}
where $SM$ is the unit normal sphere bundle to $M$ in $W$. 
The shift $s$ depends on whether the function is increasing or
decreasing and concave or convex near $x^\pm_l$ or $y^\pm_l$ 
and on the multiplicity of the orbits on the
level. Hence, the Floer homology in \eqref{eq:poz} can be
non-zero only for the range of degrees bounded by the Conley--Zehnder
indices $*$ such that $*-s=0$ and $*-s=\dim SM$. These degree ranges 
and approximate values of actions  $A(x^\pm_l)$ and $A(y^\pm_l)$ are
given in Table 1. Here we use the notations
$$
2m=\dim M\quad\text{and}\quad 2n=\codim M,
\quad\text{so that}\quad n_0=n-m+1,
$$
and the dots in the expressions for actions
stand for the terms which can be made arbitrarily small by a suitable choice
of $K_\pm$ while keeping $\max K_\pm$ and $r$ and $R$ constant; namely 
by shortening the convexity/concavity intervals. The calculation of the
actions and degree ranges from Table 1 is straightforward (but somewhat
tedious for the degrees). For the sake of completeness we
will outline it in Section \ref{sec:table} at the end of the proof.

\begin{table}
\caption{The Conley--Zehnder indices and actions for $x_l$ and $y_l$. }
 \begin{tabular}{|c|l|l|l|}

 \hline
 $\rho$ & Degrees & Actions for $K_+$ & Actions for $K_-$ \\
 \hline
 $x_1$ 
 & $[n-m+1,\, 3n+m]$
 & $ \text{max}\,K_+ + \pi R^2 \pm \ldots$ 
 & $ \text{max}\,K_- + \pi r^2 \pm \ldots$ \\
 $x_2$ 
 & $[3n-m+1,\, 5n+m]$ 
 & $\text{max}\,K_+ + 2\pi R^2 \pm \ldots$ 
 & $\text{max}\,K_- + 2\pi r^2 \pm \ldots$ \\
 \vdots 
 & \vdots 
 & \vdots
 & \vdots\\
 $x_l$ 
 & $[(2l-1)n-m+1,~ ~$
 & $\text{max}\,K_+ + l\pi R^2 \pm \ldots$ 
 & $\text{max}\,K_- + l\pi r^2 \pm \ldots$ \\
 & $~ ~ (2l+1)n+m]$ & & \\
\vdots 
 & \vdots 
 & \vdots
 & \vdots\\
 \hline
 $y_1$ 
 & $[n-m,\, 3n+m-1]$
 & $ \pi R^2 \pm \ldots $ 
 & $ \text{min}\,K_- +\pi r^2 \pm \ldots$ \\
 $y_2$ 
 & $[3n-m,\, 5n+m-1]$
 & $ 2 \pi R^2 \pm \ldots $ 
 & $ \text{min}\,K_- +2 \pi r^2 \pm \ldots$ \\
 $y_3$ 
 & $[5n-m,\, 7n+m-1]$
 & $ 3 \pi R^2 \pm \ldots $ 
 & $ \text{min}\,K_- +3 \pi r^2 \pm \ldots$ \\
\vdots 
 & \vdots 
 & \vdots
 & \vdots\\
 $y_l$ 
 & $[(2l-1)n-m,$
 & $ l\pi R^2  \pm \ldots$ 
 & $ \text{min}\,K_- + l\pi r^2 \pm \ldots $ \\
 & $(2l+1)n+m-1]$ 
 & & \\
\vdots 
 & \vdots 
 & \vdots
 & \vdots\\
\hline
\end{tabular}
\end{table}


Returning to the definition of $K_\pm$, observe that since 
$\max H > \pi R^2$, the function $K_+$ can be chosen so that 
\begin{equation}
\labell{eq:K+}
A(y^+_1)< \max H <\max K_+ < A(x^+_1) \quad\text{and}\quad  A(y^+_2)< A(x^+_1).
\end{equation}
In a similar vein, since $\min K_-\leq 0$ and $r>0$ can be taken arbitrarily
small, $K_-$ can be chosen so that
\begin{equation}
\labell{eq:K-}
A(y^-_1)< A(y^-_2) < A(y^+_1) < \max H 
\end{equation}

Now we are in a position to specify the conditions on the action
interval end-points $a$ and $b$. Namely, we only require that these points
be outside of the action spectra of $K_\pm$ and
\begin{equation}
\labell{eq:a-b}
\max K_+ < a < A(x^-_1) < A(x^+_1) < b.
\end{equation}
In particular, $b$ can be taken arbitrarily large or arbitrarily close to
$A(x^+_1)$. Thus, the actions $A(x^-_1)$ 
and $A(x^+_1)$ are necessarily in the interval $(a,\, b)$ and the interval 
automatically contains neither the points 
$$
A(y^-_1)< A(y^-_2) < A(y^+_1) < \max H  < \max K_+ ,
$$
nor $\max K_-$. The interval may contain $A(y^+_2)$, but then necessarily 
$a< A(y^+_2)< A(x^+_1)$. In addition, the interval may contain some
of the points $A(x^\pm_l)$ with $l\geq 2$ and  some
of the points $A(y^\pm_l)$ with $l\geq 3$.

\subsubsection{Showing that $\HF^{[a,\, b)}_{n_0}(K_\pm)=\Z_2$}
\labell{subsec:calc}
Let us first calculate the Floer homology for $K_+$.
By our choice of $[a,\,b)$, only the periodic orbits on
the levels 
$$
\rho=x^+_1,\, x^+_2,\,\ldots
\quad\text{and}\quad
\rho=y^+_2,\, y^+_3,\,\ldots
$$
can contribute to the homology. 

If the interval $[a,\,b)$ contains only the action $A(x^+_1)$, the
identity $\HF^{[a,\, b)}_{n_0}(K_+)=\Z_2$ follows immediately from 
\eqref{eq:poz} and the calculation
of the degrees in Table 1. (Note that $n_0$ is exactly the left
endpoint of the range of Conley--Zehnder indices for $x^\pm_1$.)

Now we argue inductively (as in \cite{cgk}) to show that the interval 
$[a,\, b)$ can be shrunk to an interval containing only $A(x^+_1)$. For 
example, let $b'$ be outside of the action spectrum of $K_+$ and such that 
$$
a< A(x^+_1)< b'<b
$$
and the interval $(b',\, b)$ contains only one point $A(x^+_l)$, $l\geq 2$, or 
$A(y^+_l)$, $l\geq 3$. (Since $A(y^+_2)< A(x^+_1)$, the action $A(y^+_2)$ 
cannot occur in the interval $(b',\, b)$.) We need to show that 
$$
\HF^{[a,\,b')}_{n_0}(K_+)\to
\HF^{[a,\,b)}_{n_0}(K_+).
$$
is an isomorphism. As can be easily seen from the table, 
$$
\HF^{[b',\,b)}_{*}(K_+)=0
\quad\text{for}\quad
\begin{cases}
*< (2l-1)n-m+1 &  \text{if $A(x^+_l)\in (b',\,b)$, $l\geq 2$},\\
*< (2l-1)n-m & \text{if $A(y^+_l)\in (b',\,b)$, $l\geq 3$}.
\end{cases}
$$
Since $n_0=n-m+1$, we have 
$\HF^{[b',\,b)}_{n_0+1}(K_+)=\HF^{[b',\,b)}_{n_0}(K_+)=0$, and the map in
question is an isomorphism by Lemma \ref{lemma:poz}. (Note that the lemma
would not apply when $n=1$ if we had $A(y^+_2)\in (b',b)$.)

Arguing inductively, we can move $b> A(x^+_1)$ to the left 
as close to $A(x^+_1)$ as we wish without changing the homology in degree
$n_0$.

Next observe that there can be some points $A(y^+_l)$, $l\geq 2$, in the
interval $(a,\,A(x^+_1))$. We repeat the same argument. Let $a'$ be
such that
$$
a< a'< A(x^+_1)
$$ 
and the interval $(a,\,a')$ contains only one of the points $A(y^+_l)$, 
$l\geq 2$. Then again using the table and applying Lemma \ref{lemma:poz}, we
see that
$$
\HF^{[a,\,b)}_{n_0}(K_+)\to
\HF^{[a',\,b)}_{n_0}(K_+)
$$
is an isomorphism. Indeed, now we need 
$n_0$ to be outside of the range of degrees for $y_l^+$, i.e.,
$n_0< (2l-1)n-m$ for $l\geq 2$, which is clearly
true with $y^+_2$ being the worst case scenario. (Then $n_0-1$ is 
automatically outside of the range of degrees.)
Therefore, we can move $a< A(x^+_1)$ to the
right as close to $A(x^+_1)$ as we wish, without changing the homology
in degree $n_0$.

Hence, for the above choice of $[a,\,b)$, the homology
$\HF^{[a,\,b)}_{n_0}(K_+)$ is the same as when $[a,\,b)$ contains
only $A(x^+_1)$ and we conclude that this group is $\Z_2$.

Regarding the calculation of $\HF^{[a,\, b)}_{n_0}(K_-)$, note that 
the points $y^-_l$ and $y^+_l$ have the same ranges of Conley--Zehnder 
indices. This is also true for the points $x^-_l$ and $x^+_l$. Furthermore,
now not only $A(y^-_1)$ but also $A(y^-_2)$ is outside of the 
interval $[a,\,b)$.
Hence, the argument we have used for $K_+$ translates word-for-word to 
the calculation of the Floer homology for $K_-$ and, therefore,
$\HF^{[a,\, b)}_{n_0}(K_-)=\Z_2$. 

\subsubsection{Monotone homotopy and the isomorphism  
$\HF^{[a,\, b)}_{n_0}(K_+)\to \HF^{[a,\, b)}_{n_0}(K_-)$}
Before describing the monotone homotopy, let us observe that without
loss of generality the functions $K_\pm$ can be assumed to have approximately
equal slopes and hence equal number of periodic levels $x^\pm_l$ and $y^\pm_l$.
This can be achieved by either starting with functions satisfying this
requirement or by increasing the slope of one of them through a monotone
homotopy. (In the latter case, new periodic levels are created with
actions inside of $(a,\,b)$, but the homotopy can be arranged so that $a$ and
$b$ stay away from the action spectrum.)

The monotone homotopy $K_s$, $s\in [0,1]$, from $K_0=K_+$ to $K_1=K_-$ is
shown in Fig.\ 2. In the course of this homotopy, the bottom part of $K_+$ 
moves down eventually reaching $\min K_-$. Then, at the second stage of 
the homotopy, $\max K_+$ moves down to $\max K_-$ in a monotone fashion and 
the linear part of the function moves to the left. It is easy to see
that the homotopy can be arranged so that the 
periodic energy levels considered above persist under the homotopy
and remain Morse--Bott non-degenerate. Thus, $x^+_l$ (or 
$y^+_l$) moves to $x^-_l$ ($y^-_l$, respectively) through a family 
of periodic levels $x^s_l$ ($y^s_l$, respectively) and $x^s_l$ and $y^s_l$
are smooth functions of $s\in [0,1]$. Furthermore, as is clear from Fig. 2,
the actions $A(x^s_l)$ and $A(y^s_l)$ can be assumed to be monotone 
decreasing functions of $s$.

\begin{figure}[htbp]
\caption{The homotopy from $K_+$ to $K_-$}
\begin{center}
\input{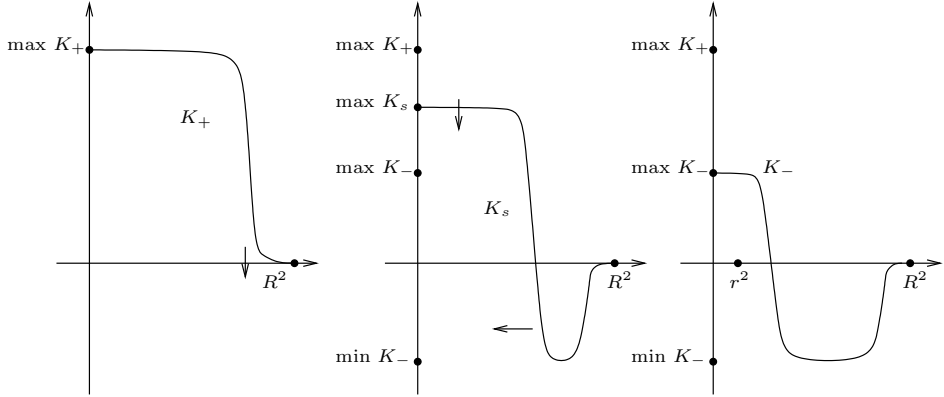}
\label{figure:homotopy}
\end{center}
\end{figure}

Our goal is to show that $\HF^{[a,\,b)}_{n_0}(K_s)$ remains constant during
the homotopy. If $a$ and $b$ were always away from the action spectrum 
of $K_s$, this would follow (for all degrees) from the homotopy invariance 
of Floer homology. However, this is not the case and we
have to analyze the behavior of periodic orbits under the homotopy more
closely.

For the sake of simplicity, let us first assume that 
$b$ is to the right of the action spectra of all $K_s$. 
(Recall that $b$ can be taken arbitrarily large.) Then
the actions $A(x^s_l)$, $l\geq 1$, are inside the interval $(a,\,b)$
for all $s\in [0,1]$.
Furthermore, $A(y_1^s)$ is outside of $(a,\, b)$ for all $s$.
Thus we have
$$
A(y^s_1) < a < A(x^s_1)< A(x^s_2)< A(x^s_3)< \ldots < b 
$$
for all $s$.
The actions $A(y_l^s)$ with $l\geq 2$ may cross the left end-point of 
the interval. Finally, new periodic levels are necessarily created within the 
interval where $K_s$ is increasing. However, the orbits on these 
levels have negative actions and hence do not contribute to the homology.

Hence, we only need to examine the effect of $A(y^s_l)$, $l\geq 2$, crossing
the left end-point $a$ of the interval at some moment $s_0$. This effect 
is the same as when moving $a$ through the action value $A(x_l^{s_0})$ 
in the opposite direction. The exact sequence argument used in the
previous section applies in this case, and the homology in degree
$n_0$ remains unchanged.  

For the sake of completeness, let us outline a rigorous proof of this fact.
Let $A(y^{s_0}_l)=a$. Since $A(y^s_l)$ is a monotone decreasing 
function of $s$, without loss of generality we may assume that there exists
a small interval $I=(s_1,\, s_2)$ containing $s_0$ and a small interval
$(a_1,\, a_2)$ containing $a$ such that for every $s\in (s_1,\,s_2)$ 
\begin{itemize}

\item $a_1 < A(x^{s_1}_l) < a=A(x^{s_0}_l) < A(x^{s_2}_l) <a_2$ and

\item $A(y^s_l)$ is
the only point of the action spectrum of $K_s$ in $(a_1,\,a_2)$.

\end{itemize}
Then the monotone homotopy map 
$$
\HF_{n_0}^{[a,\,b)}(K_{s_2}) \to \HF_{n_0}^{[a,\,b)}(K_{s_1})
$$ 
factors as
$$
\HF_{n_0}^{[a,\,b)}(K_{s_2}) 
\to\HF_{n_0}^{[a_1,\,b)}(K_{s_2}) 
\to\HF_{n_0}^{[a_1,\,b)}(K_{s_1}) 
\to \HF_{n_0}^{[a,\,b)}(K_{s_1}).
$$
Here, the first map is an isomorphism because $K_{s_2}$ does not have
action spectrum values in $[a_1,\,a]$. The second map is an isomorphism
by homotopy invariance. Finally, the third map is an isomorphism by the exact
sequence argument. Namely, first note that $y^s_l$ has the same range of
indices as $y^\pm_l$. Then, using Table 1, we see that
$\HF^{[a_1,\,a)}_{n_0-1}(K_{s_1})=\HF^{[a_1,\,a)}_{n_0}(K_{s_1})=0$ and
hence, by Lemma \ref{lemma:poz},  the third map is also an isomorphism. 

A similar argument shows that the monotone homotopy induces an isomorphism
of the Floer homology in degree $n_0$ for any $b> A(x^+_1)$.

This concludes the proof of Proposition \ref{prop:floer} and
of Theorems \ref{thm:cap-nbd} and \ref{thm:main}. 
In the remaining section of the proof, we outline 
the calculation of the actions and degree ranges
given in Table 1 

\subsubsection{Calculation of the actions and degree ranges from Table 1.}
\labell{sec:table}
The actions $A(x^\pm_l)$ and $A(y^\pm_l)$ are easy to determine. For
instance, the periodic orbits of $K_+$ on the level $x^+_l$ are the
$l$-iterated Hopf circles in the fibers of the symplectic normal
bundle to $M$.  These circles have radius $\sqrt{x^+_l}$, i.e.,
approximately $R$, and are traversed in the negative
direction, for $K_+$ is decreasing.  Thus, the symplectic area bounded
by the circles is approximately $-l\pi R^2$ and the value of $K_+$
is near $\max K_+$. Hence, the action is approximately $\max
K_+ + l\pi R^2$.

Let us now turn to the ranges of Conley--Zehnder indices.
For the sake of brevity we indicate only the main steps of
the calculation. Assume first that $M$ is the origin in $W=\C^n$.
Considering an explicit perturbation of $K_\pm$ near $x_1^\pm$ and
$y_1^\pm$ and using the definition of the Conley--Zehnder index, it
is not hard to see that the range of degrees is $[-n+1,\,n]$ when
$K_\pm$ is convex near the level (i.e., for $y_1^\pm$) and 
$[-n,\,n-1]$  when $K_\pm$ is concave near the level (i.e., for $x_1^\pm$).
These perturbations can be obtained by taking a small non-degenerate
quadratic form on $\C P^{n-1}$ and then making a time-dependent perturbation
supported near now-isolated periodic orbits. The indices are easy to
calculate in the trivializations arising from those of the tangent spaces
to $\C P^{n-1}$ at the projections of the orbits. Then, passing 
to the standard trivialization of $\C^n$ amounts to shifting
the range of indices by $2n$. 
(We refer the reader to, e.g., \cite{sa} for the definition of the
Conley--Zehnder index, its properties, and further references. The reader
should keep in mind that the sign convention of \cite{sa} is different
from the one used here and this difference affects Conley--Zehnder indices.
For example, for the action defined by \eqref{eq:action}, a non-degenerate 
critical point with Hessian $S$ of an autonomous 
$C^2$-small Hamiltonian has Conley--Zehnder index $-\text{signature}(S)/2$,
but not $\text{signature}(S)/2$ as in \cite{sa}.)
For iterated orbits, i.e., the levels $x_l^\pm$ and $y_l^\pm$, these
ranges are further shifted by $2n(l-1)$. This again can be seen as a result 
of a trivialization change. Finally, to deal with the general case of a tubular
neighborhood of $M\subset W$, we perturb $K_\pm$ by adding to it a 
$C^2$-small Morse function on $M$. This results in addition of the
interval $[-m,\, m]$ to the range of Conley--Zehnder indices.

\subsection{Proof of Theorem \ref{thm:time-dependent2}} 
As in the proof of Theorem \ref{thm:main}, first observe that it is
sufficient to prove the theorem for a tubular neighborhood $U_R$ of $M$.
Reasoning as above, we will find functions $K_+$ and $K_-$ (independent of 
time), supported in 
$U_R$, such that $K_-<H<K_+$ and 
\begin{equation}
\labell{eq:t-dependent}
\HF^{[a,\,b)}(K_+)\to \HF^{[a,\,b)}(K_-)
\end{equation}
is an isomorphism for some constants $0<a<b$, not in the action
spectra of $K_\pm$, with $a>0$ being arbitrarily small. Then the theorem 
will follow from \eqref{eq:t-dependent}.

The function $K_+$ has the same shape as its counterpart in Fig.\ 1 with
the only modification that at the top part $K_+$ is now slowly monotone
decreasing with constant slope. Thus $K_+$ has a Morse--Bott non-degenerate
maximum along $M$ (with small eigenvalues). It is clear that such a function 
can be chosen so that
$H\leq K_+$. The function $K_+$ has one-periodic levels at $\rho$ equal to
$x_1,\ldots$ and $y_1,\ldots$ with actions and Conley--Zehnder indices
as in Table 1. In addition to this, $M$ is a Morse--Bott non-degenerate
set of one-periodic
orbits of $K_+$. The range of indices for $M$ is $[n-m,n+m]$ and the action
is $\max K_+$. As is clear from Table 1, the action spectrum
of $K_+$ is strictly positive except the action on the trivial orbits 
with $K_+=0$.

Let $K_-(\rho)=\eps K_+(\eps \rho)$, where $0<\eps<1$. In other words,
the graph of $K_-$ is obtained from the graph of $K_+$ by scaling by $\eps$.
Since $H\geq 0$ and $\min_{S^1\times M} H>0$, we clearly have $K_-\leq H$ 
if $\eps>0$ is small enough.

Consider now the homotopy $K_s(\rho)=s K_+(s \rho)$ with $s\in [\eps,\,1]$ 
from $K_+$ to $K_-$. Again, all action values of $K_s$ (on periodic
levels where $K_s> 0$) are separated from zero. It follows from homotopy
invariance of Floer homology that for any sufficiently large $b$ and any
sufficiently small $a>0$, the monotonicity map
$$
 \HF^{[a,\,b)}(K_+)\to \HF^{[a,\,b)}(K_-)
$$
is an isomorphism in all degrees. 

Now it is sufficient to show that $\HF^{[a,\,b)}(K_-)\cong H_*(M;\Z_2)$.
To this end, consider
the monotone (increasing) homotopy $F_s(\rho)=\eps K_+(s\rho)$ with 
$s\in [\eps,\, r]$. This homotopy  begins with $F_\eps=K_-$ and ends with
$F_r\geq K_-$. (Note that $F_s$ is obtained by dilating the graph of $K_-$
along the $\rho$-axis; the homotopy is defined only if $r>\eps$ is small 
enough, e.g., $r<R$.) In the course of the homotopy $F_s$ from
$s=\eps$ to $s=r$, the periodic levels $x_l$ and $y_l$ come close to each 
other, collide, and disappear. No new periodic levels with action close to 
zero are created.
Thus
$$
 \HF^{[a,\,b)}(F_r)\to \HF^{[a,\,b)}(K_-),
$$ 
is an isomorphism
if $b$ is sufficiently large and $a>0$ is sufficiently small.
If $\eps>0$ is small enough and, say, $r=R/2$, all periodic levels
$x_l$ and $y_l$ of $F_r$ are destroyed. In other words, $F_r$ does not have 
one-periodic orbits other than the critical manifold $M$ (in the region 
where $F_r>0$). By Po\'zniak's theorem,
$$
 \HF^{[a,\,b)}(F_r)\cong H_*(M;\Z_2)
$$
with some shift in degrees, which completes the proof of the theorem.

\begin{Remark}
\labell{rmk:$H>0$}
It should be possible to show directly that 
$\HF^{[a,\,b)}(K_\pm)\cong H_*(M;\Z_2)$. Note,  however, that
since the ranges of Conley--Zehnder indices of $x_1$ and $y_1$ are close
to that for $M$, the exact sequence argument utilized in the proof of
Proposition \ref{prop:floer} does not apply. Hence, a more delicate
reasoning is required, e.g., a calculation of the differential in the Floer
complex (cf. \cite{fhw}). Although this approach is more involved than the 
homotopy argument above, it can perhaps be utilized to relax the hypotheses 
of Theorem \ref{thm:time-dependent}.
\end{Remark}

\subsection{Discussion}
As has been pointed out, Theorem \ref{thm:main} does not hold if the 
assumption that $H$ is constant near $M$ is replaced by the weaker assumption 
that $H$ is constant on $M$. Namely, as the example below shows,
for any $R>0$ and $c>0$ there exists
a smooth non-negative function $H$, constant on $M$ and supported in $U_R$, 
such that $c=\max H$
and all non-trivial one-periodic orbits of $H$ have action less than or 
equal to $c$. The Hamiltonian flow
of $H$ has numerous non-trivial one-periodic orbits (it should, by Theorem
\ref{thm:min}, have a fast non-trivial periodic orbit), but these orbits 
are not detected
by the Floer homology with the range of actions greater than $\max H$.

\begin{Example}
\labell{exam:actions}
Fix $b=\pi k+\pi/2$, where $k$ is a positive integer, large enough so that 
$c/b$ 
(the solution of $c-by=0$) is in the interval $[0,R^2]$. Let $H=f(\rho)$, 
where $f$ is a smooth function of the form 
$$
f(y)=\begin{cases}
c-by &\text{for $y\in [0,\delta_-]$,}\\
\text{monotone decreasing} & \text{for $y\in (\delta_-,\delta_+)$,}\\
0 & \text{for $y\in [\delta_+,\infty)$.}
\end{cases}
$$
Here, $[\delta_-,\delta_+]$ is an arbitrarily small interval containing 
$c/b$ and contained in $(0,R^2)$.

Non-trivial one-periodic orbits of $H$ occur on the sphere bundles $\rho=y_l$, 
where $y_l$ 
is the solution of the equation
$$
f'(y_l)=-\pi l, \quad l\in \Z_+.  
$$
Since $b$ is not an integer multiple of $\pi$, we have $y_l\in 
(\delta_-,\delta_+)$ and
$l\in [1,b/\pi)$. All orbits on the level $\rho=y_l$ have action 
\begin{eqnarray*}
A_l &=& \pi l y_l + f(y_l)\\
&=& \pi l y_l +(c-by_l) + \big(f(y_l)-(c-by_l)\big)\\
&=& c+ (\pi l-b)y_l + \big(f(y_l)-(c-by_l)\big).
\end{eqnarray*}
It is easy to see that $(\pi l-b)y_l< -\pi \delta_-/2$ and 
$f(y_l)-(c-by_l) < c-b\delta_+$ since $y_l\in (\delta_-,\delta_+)$. Thus,
$$
A_l-c<  -\pi \delta_-/2+ ( c-b\delta_+) <0
$$
if $\delta_+$ is taken sufficiently close to $c/b$ while $\delta_-$ is fixed. 
Therefore,
$A_l< c$ for all $l$.

Note also that when $c>0$ is small, all one-periodic orbits of $H$
are trivial.
\end{Example}

\begin{Example}[Polterovich]
\labell{exam:polt}
In this example, due to Leonid Polterovich, we show that 
in contrast with Theorem \ref{thm:main}, the action value in 
Theorem \ref{thm:time-dependent} cannot be bounded from below via 
$\min_{S^1\times M} H$ even when $M$ is a point.
To be more specific, \emph{for any $\eps>0$ and $C>0$, there exists a 
non-negative Hamiltonian $H\colon S^1\times \R^{2n}\to \R$
with arbitrarily small support such that $\min H|_{S^1\times 0}\geq C$ and
every one-periodic orbit of $H$ has action less than $\eps$}.

Indeed, first
note that the condition that $H$ is periodic in time can be dropped by
\cite[Proposition 2.1.3]{bps}. (Here some extra care is needed to keep
$H$ non-negative.) Let $K$ be a $C^2$-small non-negative
Hamiltonian supported in some ball, such that the time-one flow of $K$
displaces a ball $B$ centered at the origin and such that
$\max K<\eps$. Let $F$ be a bump-function supported in $B$ with $F(0)\geq C$.
Consider the compactly supported Hamiltonian $H_t$ generating the
time-dependent flow $\psi_t\varphi_t$, where $\psi_t$ is the flow of $F$ and
$\varphi_t$ is the flow of $K$. It is easy to see that $H_t(0)\geq C$ and
the only one-periodic orbits of $H_t$ are the critical points of $K$.
Hence, all non-zero action values for $H_t$ are less than $\eps$.

A similar construction with $K$ now being non-positive shows that
\emph{for any $\eps>0$ and $C>0$, there exists a  Hamiltonian 
$H\colon S^1\times \R^{2n}\to \R$
with arbitrarily small support such that $\min H|_{S^1\times 0}\geq C$,
every one-periodic orbit of $H$ has a non-positive action, and the 
time-one flow of $H$ is non-trivial}.
\end{Example}

\begin{Remark}
We conclude this discussion by pointing out the location of
orbits from Theorems \ref{thm:min} and \ref{thm:time-dependent}
for the function $H$ from Example \ref{exam:actions}.
The orbits from Theorem \ref{thm:time-dependent} lie on the critical manifold
$M$. The orbits from Theorem \ref{thm:min}, 
obtained via Theorem \ref{thm:main}, lie on the level set $\rho=x_1$, in the
notations of Table 1. Finally,
the orbits which make the homological capacity from \cite{cgk} non-vanish
are located on the level $\rho=y_1$ (see also Section \ref{sec:other-cap}).
A suitably adapted version of the Schwarz action selector (see \cite{FS,sc})
will pick either $\max K_+$ or $A(y_1^+)$, whichever is smaller.
\end{Remark}

\section{Capacity: proofs and remarks}
\labell{sec:capacity-pf}

\subsection{Proof of Theorems \ref{thm:min} and \ref{thm:HZ-modify}} 
The key to the proof of these theorems is the following elementary observation,
essentially contained already in \cite[p. 184]{hz:book},
which allows one to cut off an autonomous Hamiltonian without creating new
fast periodic orbits.

\begin{Lemma}
\labell{lemma:cutoff}
Let, in the notations of Section \ref{sec:defHZ}, $H\in \tilde{\CH}_0(V,Z)$ 
be such that the flow of $H$ has no (contractible) non-trivial fast periodic 
orbits. Then, for any $C$ with $0<C < \min_Z H$, there exists a function 
$K\in \CH(V,Z)$ such that $C \leq \max K < \min_Z H$ and the flow of $K$ has 
no non-trivial fast periodic orbits.
\end{Lemma}

\begin{Remark}
This lemma and Lemma \ref{lemma:positive} stated below still hold when
\emph{all} fast non-trivial periodic orbits of $H$ and $K$ are replaced by 
\emph{contractible in $V$} non-trivial fast periodic orbits.
\end{Remark}

Theorem \ref{thm:min} immediately follows from Lemma \ref{lemma:cutoff}
and Theorem \ref{thm:main}; Remark \ref{rmk:period} is a consequence
of Lemma \ref{lemma:cutoff}, Theorem \ref{thm:main}, and the definition
of relative Hofer--Zehnder capacity. 

For the sake of completeness, we outline the proof of Lemma \ref{lemma:cutoff}.

\begin{proof} Without loss of generality, we may assume that $H\neq 0$.
Pick a constant $C<\min_Z H$. Let $\eps>0$ be sufficiently small 
so that  $0< C-\eps$ and $C+\eps< \min_Z H$.
Then the function $K$ is obtained by cutting $H$ off at the level $H=C$ and
then smoothening up the resulting function. More precisely, let 
$f\colon [C-\eps,C+\eps]\to \R$ be a function with $0\leq f'\leq 1$ such that
$f(y)=y$ for $y$ near the left end-point of the interval and 
$f(y)=C$ for $y$ near the right end-point of the interval. Set
$$
K(z)=
\begin{cases}
H(z) & \quad\text{when $H(z)\leq C-\eps$,}\\
f(H(z)) & \quad\text{when $C-\eps\leq H(z)\leq C+\eps$,}\\
C & \quad\text{when $C+\eps\leq H(z)$.}
\end{cases}
$$
Then $K$ is constant (and equal to $C$) near $Z$ and equal to $H$ near
the boundary of $\supp H$, i.e., $K\in\CH(V,Z)$. Furthermore, $\max K=C$ and 
$K$ does not have fast periodic orbits since $|f'|\leq 1$.
\end{proof}

Let us prove the first assertion of Theorem \ref{thm:HZ-modify}, i.e.,
$$
\CHZ(V,Z)=\inf E,
$$
where
$$
E=\{C>0\mid\text{every $H\in \TCH_C(V,Z)$ 
has a non-trivial fast periodic orbit}\}.
$$ 
Note that $E$ is a semi-infinite interval  of the form
$[\inf E,\infty)$ or $(\inf E,\infty)$ or the empty set. 
With this in mind, we have
$\CHZ(V,Z)\leq \inf E$, by definition. (Indeed, if $H\in \CH(V,Z)$ is
such that all periodic orbits of $H$ are slow, then $\max H$ must be
outside of this interval.) To prove the opposite inequality, assume that
$H\in \TCH_C(V,Z)$ has no fast non-trivial periodic orbits. Then $C$ is
in the complement of $E$. By
Lemma \ref{lemma:cutoff}, for any small $\eps>0$
there exists $K\in \CH(V,Z)$ without fast non-trivial periodic orbits and
such that $\max K> C-\eps$. Thus, $\CHZ(V,Z)\geq C$ and, since the complement
of $E$ is also an interval, $\CHZ(V,Z)\geq \inf E$.

Likewise, the second assertion of Theorem \ref{thm:HZ-modify}
is an easy consequence of

\begin{Lemma}
\labell{lemma:positive}
Let $H\in \CH(V,Z)$ ($H\in \TCH_C(V,Z)$ with $C>0$)
be such that the flow of $H$ has
no non-trivial fast periodic orbits.
Then, there exists a non-negative function $K\in\CH(V,Z)$
(respectively, $K\in \TCH_C(V,Z)$)
with $\max K = \max H$ having no non-trivial fast periodic orbits. 
\end{Lemma}

\begin{proof}
Recall that a $C^2$-small function (with support in a given compact set)
does not have fast non-trivial periodic orbits; see, e.g., 
\cite[pp.\ 185, 200]{hz:book}. Using this fact, it is not hard to modify
$H$ so that it becomes non-negative near the boundary (in the sense of
general topology) of $\supp H$, without
changing $\max H$ or creating fast periodic orbits. We denote the resulting
function by $H$ again. 

Let $c>0$ be a small regular value of $H$. Denote by $Y_l,\ldots, Y_k$
the connected components of $\{H\leq c\}$ on which $H$ assumes negative
values. Let also $\eps>0$ be so small that $c-\eps>0$ 
and different connected components $Y_l,\ldots, Y_k$ are contained in different
connected components of $\{H\leq c+\eps\}$. (As a consequence, 
$c+\eps<\max H$.) Now, separately for each $Y_i$, we cut off $H$ along 
$\partial Y_i$ and smoothen it up exactly as in the proof of Lemma 
\ref{lemma:cutoff}. It is clear that the resulting function $K$ has the
required properties.
\end{proof}

This concludes the proof of the theorems.

\subsection{Almost existence theorem}
\labell{sec:alm-ext}
The proof of the almost existence theorem, both the standard 
(see \cite{hz:book}) and relative (Theorem \ref{thm:ae})
versions, can be broken down into the following two results which may be of 
independent interest.

\begin{Proposition}
\labell{prop:ineq}
Let $U$ be an open subset of a symplectic manifold $V$ with compact closure 
$\bar U$ and let $Z\subset U$ be compact. Then
$$
\CHZ(V,\bar U)\leq \CHZ(V,Z)-\CHZ(U,Z).
$$
In particular,
$$
\CHZ(V,\bar U)\leq \CHZ(V)-\CHZ(U).
$$
\end{Proposition}

\begin{Example}
Let $B_R$ be the ball of radius $R>0$ in $\R^{2n}$, centered at the
origin. An argument similar to the proof of Theorem \ref{thm:main}
shows that $\CHZ(B(R),\bar{B}(r))=\pi (R^2-r^2)$. It is not clear whether
the inequalities in Proposition \ref{prop:ineq} are not in fact 
equalities.
\end{Example}

To further deal with the almost existence problem, let us, following
\cite{hz:book}, analyze the existence of periodic orbits on a given
hypersurface $S$ bounding a domain $U$ in $V$.
Let $S_\eps$ be a thickening of $S=S_0$ in $V$. Denote by
$U_\eps$ the domain bounded by $S_\eps$. We say that $S$ has
\emph{relative Lipschitz type} if
$$
\limsup_{\eps \searrow 0} \frac{\CHZ(U_\eps,\bar U)}{\eps} <\infty.
$$
It is easy to see that this is a well-defined property, i.e.,
independent of the choice of the thickening $S_\eps$. Note that, by
Proposition \ref{prop:ineq}, a Lipschitz type hypersurface (see 
\cite{hz:book}) is automatically of relative Lipschitz type.

\begin{Proposition}
\labell{prop:Lipschitz}
A hypersurface of relative Lipschitz type carries a closed characteristic.
\end{Proposition}

We leave both of these propositions without proofs, for the arguments are 
implicitly contained in \cite[Section 4.2]{hz:book}. Similar results
hold for the restricted Hofer--Zehnder capacity $\CW$.

\begin{Remark}
It is not clear if there exist hypersurfaces of 
relative Lipschitz type which are not of Lipschitz type.
\end{Remark}

\subsection{Concluding remarks}
\labell{sec:other-cap}
Arguments used to prove finiteness of a capacity or the
existence of periodic orbits can often be turned
into a definition of a new capacity, bounding the original one from
above; see, e.g., \cite{bps,cgk,fhw}. The proof of Theorem \ref{thm:main}
is no exception. This leads to the notion of a restricted homological 
capacity $\wsh$ which bounds $\CW$ from above and is equal to $\pi r^2$
on $(U_r,M)$. (Here we assume that the ambient manifold $W$ is geometrically
bounded and symplectically aspherical; $M$ is a closed symplectic
submanifold of $W$.) The definition of this homological capacity is
a straightforward (but cumbersome) axiomatization of the proof, and we omit
it here. The only advantage the capacity $\wsh$ seems to have over
$\CW$ is that $\wsh$ gives some information about the actions of periodic
orbits. It is not clear, however, how to use this extra information.

Let us briefly discuss the relation of the capacity $\CHZ$ or $\CW$ with
the capacities introduced in \cite{bps} and \cite{cgk}. 

The restricted relative
capacity of \cite{cgk} is a relative version of the homological
capacity from \cite{fhw}. The finiteness of this capacity results
in the ``nearby existence theorem'' -- the existence of periodic
orbits on a dense set of levels (see \cite{cgk}) -- but falls short 
of leading to the almost existence theorem. We are not aware of any 
inequalities relating this capacity with $\CW$ or $\wsh$.

For the trivial homotopy class, the relative capacity (homological
or not) of \cite{bps} does not allow one to control whether the periodic orbit
detected by the capacity is trivial or not (cf. Example \ref{exam:actions}).
Thus, for the trivial homotopy class this capacity appears to be unrelated
to $\CHZ$. On the other hand, when the homotopy class is non-zero, the orbit
is automatically non-trivial and in this case the capacity of \cite{bps}
gives an upper bound for $\CHZ$ (but not $\CW$).

\begin{Example}
\labell{exam:Lagrangian}
Let $Z$ be a closed Lagrangian submanifold of a geometrically bounded
symplectically aspherical manifold $W$. When $Z$ is a torus or admits a metric
of negative sectional curvature, $\CHZ(U,Z) <\infty$, where $U$ is a small
neighborhood of $Z$, as immediately follows from the results of \cite{bps}.
The results of \cite{hv-ct,vi} also suggest, but apparently do not imply,
that $\CHZ(U,Z)<\infty$ in general for a closed Lagrangian submanifold.

Note also that in this case $\CW(U,Z)=\infty$, when $Z$ admits a metric without
contractible geodesics and $\pi_1(Z)\to \pi_1(W)$ is a monomorphism.
The restricted capacity $\CW$ is specifically ``tuned up'' to detect
contractible periodic orbits near a compact submanifold and
Theorem \ref{thm:main} guarantees the existence of such orbits near
a symplectic submanifold. For a general submanifold, such periodic orbits 
may fail to exist, and the restricted capacity $\CW$ may be infinite even when 
non-contractible orbits exist in abundance.
\end{Example}

As was pointed out by L. Polterovich, one may expect that 
$\CHZ(V,M)=\CHZ(V)$ when $M$ is a symplectic submanifold or even
when $\omega|_M\neq 0$ and some natural topological conditions
hold, e.g., the normal 
bundle to $M$ in $V$ admits a non-vanishing section; cf. \cite{pol1}
and Example \ref{exam:CHZ-c0}.
(Note that finiteness of $\CHZ(V)$ or $\CW(V,\point)$ leads to
a non-relative, stronger than Theorem \ref{thm:ae}, almost existence 
theorem in $V$.) Along these lines, E. Kerman has recently shown,
\cite{Ke:new}, that
$\CW(U_r,M)=\CW(U_r)$ for $r>0$ small, provided that $M$ is a closed rational
symplectic submanifold of $W$ and the homology of the unit normal bundle 
to $M$ splits.

\begin{Remark}[Capacity of the cylinder and ellipsoids]
It is not surprising that the method used in the proof of
Theorem \ref{thm:main} readily lends itself for calculations of
capacities of some other manifolds. For example,
one can easily recover the well-known calculation of the Hofer--Zehnder 
capacity of ellipsoids (see, e.g., \cite{hz:book}),
bypassing the calculation for the cylinder, and then
derive from it the result for the cylinder. Namely, consider the solid 
ellipsoid 
$U=\{z\in \C^n\mid Q(z)<1\}$, where
$$
Q(z)=\frac{|z_1|^2}{r_1^2}+\cdots+\frac{|z_n|^2}{r_n^2}
$$
and $0<r_1\leq\ldots\leq r_n$. As is easy to see, $\CHZ(U)\geq \pi r_1^2$.
On the other hand, arguing as in the proof of Proposition \ref{prop:floer} 
with $K_\pm$ being now functions of $Q$, one can show that 
$\CHZ(U)\leq \pi r_1^2$ and, hence, $\CHZ(U)= \pi r_1^2$. 

Note now that $\CHZ(V)=\sup_U \CHZ(U)$, where the supremum is taken
over all $U\subset V$ with compact closure. (The same is true for relative
capacities, as long as $Z\subset U$.) Applying this to an exhaustion
of a symplectic cylinder by ellipsoids, we obtain that
$\CHZ(B^{2n}_r\times \R^{2m})=\pi r^2$, where $B^{2n}_r$ is the ball
of radius $r>0$.
\end{Remark}

\section{Constructions of counterexamples}
\labell{sec:counterexamples-prfs}

In this section we prove  Propositions 
\ref{prop:non-exist-R2n}, \ref{prop:non-exist}, and \ref{prop:non-exist-CPn}.
Before getting into technical details of the proofs let us outline
the basic line of reasoning in these constructions. All three arguments 
are similar concatenations  of the following standard steps:
\begin{itemize}
\item Finding or creating, if not readily available, a function with a finite
number of periodic orbits on a given level.

\item Eliminating periodic orbits on this level.

\item Applying the preceding two steps to a sequence of levels converging 
to zero.

\item Smoothening up the resulting function at zero.
\end{itemize}
Hence, we describe the proofs with a varying degree of detail, focusing
only on essential points.

\subsection{Smoothing lemma}
The last step is identical in all three proofs and we state it here as 
a lemma.

\begin{Lemma}
\labell{lemma:smoothing}
Let $Z$ be a closed submanifold of a manifold $W$ and let 
$F\colon  W\to \R$ be a function such that
\begin{enumerate}
\item $F$ is non-negative, continuous, and vanishes on $Z$;
\item $F$ is $C^m$-smooth on the complement of $Z$ for some 
$0\leq m\leq \infty$.
\end{enumerate}
Then there exists a monotone
increasing $C^\infty$-smooth function 
$\varphi\colon [0,\infty)\to [0,\infty)$ with 
$\varphi(0)=0$ and $\varphi'(y)> 0$ for $y>0$, and 
such that $H=\varphi\circ F$ is $C^m$-smooth.
\end{Lemma}

\begin{Remark}
We will need this lemma only in the cases where $m=2$ and $m=\infty$.
Below we prove it for $m=\infty$. The $m=2$ case requires only obvious
modifications.
\end{Remark}

\begin{proof}\footnote{The authors are
grateful to Ant\'ony Serra for the proof of the lemma; \cite{se}.} 
Fix a compact neighborhood $K$ of $Z$.
Let $\varphi_l\colon [0,1] \to [0,\infty)$, $l=0,1,2,\ldots$, be a sequence 
of smooth function such that
\begin{enumerate}
\item[(S1)] 
$\varphi_l(y)=0$ when $0\leq y\leq b^-_l$ and 
$\varphi_l(y)=c_l$ when  $b^+_l\leq y \leq1$,
\item[(S2)]  
$\varphi'_l(y)>0$ when $b^-_l< y< b^+_l$
\end{enumerate}
for some intervals  $(b^-_l,b^+_l)\subset (0,1)$ and constants $c_l>0$,
chosen so that the adjacent intervals overlap, $b^+_l\to 0$, and
\begin{equation}
\labell{eq:varphi}
\parallel\varphi_l\parallel_{C^{l}([0,1])}\leq \frac{1}{2^l}
\quad\text{and}\quad
\parallel\varphi_l\circ F \parallel_{C^{l}(K)}\leq \frac{1}{2^l}.
\end{equation}
The first condition in \eqref{eq:varphi} guarantees that
$\varphi=\sum\varphi_l$ is a smooth function on a neighborhood of zero,
vanishing at zero,  and such that
$\varphi'(y)>0$ for $y>0$, by (S2) and since the intervals overlap.
Let us extend this function from a neighborhood of zero to $[0,\infty)$
so that the resulting function $[0,\infty)\to [0,\infty)$,
denoted again by $\varphi$, is smooth and still has these properties.
Then, by the second condition of 
\eqref{eq:varphi}, $\varphi\circ F$ is smooth, i.e., $\varphi$ is the required
function. 
\end{proof}

\subsection{Proof of Proposition \ref{prop:non-exist-R2n}}
We start with an irrational positive-definite quadratic 
form $G\colon \R^{2n}\to \R$.
Clearly, every level $G=a$ has exactly $n$ distinct periodic orbits.

Assume first that $n>2$ so that $\dim \R^{2n} >4$. By inserting symplectic 
plugs as in 
\cite{gi:seifert95,gi:seifert97,gi:bayarea,herman,ke2}, we can eliminate
periodic orbits on a sequence of levels $ G=a_k$ with $a_k\to 0$. As a
result, we obtain a function $F\colon  \R^{2n}\to \R$ with the following
properties:
\begin{itemize}
\item $F$ meets the requirements of Lemma \ref{lemma:smoothing}
with $W=\R^{2n}$, $Z=\{0\}$, and $m=\infty$;
\item $F$ attains its absolute minimum at zero;
\item the levels $F=a_k$ carry no periodic orbits of the flow of $F$.
\end{itemize}
By Lemma \ref{lemma:smoothing}, there exists a monotone
increasing smooth function $\varphi\colon [0,\infty)\to [0,\infty)$ with 
$\varphi(0)=0$ and $\varphi'(y)> 0$ for $y>0$ and 
such that $H=\varphi\circ F$ is $C^\infty$-smooth. Then it remains to
set $\eps_k=\varphi(a_k)$.

When $n=2$, the argument is similar, but the result of \cite{gg1,gg2} is
applied to eliminate periodic orbits. In this way we obtain a function $F$ 
which is only
$C^2$-smooth on the complement of the zero section. As a consequence, the 
function $H$ is also only $C^2$-smooth. This completes the proof
of Proposition \ref{prop:non-exist-R2n}.

\subsection{Proof of Proposition \ref{prop:non-exist}} 
Here, we start with a non-symmetric Finsler metric $G\colon T^*S^n\to \R$ 
whose geodesic flow has a finite number of closed geodesics.
A metric with this property has been constructed by A. Katok, \cite{ka};
see also \cite{zi}. The rest of the proof is identical to that of
Proposition \ref{prop:non-exist-R2n}. 

To deal with the case of the twisted geodesic flow for $\sigma\neq 0$,
we observe that the flow of $G$ can be viewed as the twisted geodesic
flow of the standard metric on $S^n$ for some exact non-zero magnetic
field $\sigma$. The construction is finished in the same way as for 
the geodesic flow of $G$. 

\subsection{Proof of Proposition \ref{prop:non-exist-CPn}}
The starting point of the construction is again a function $G$ whose
flow has only finitely many periodic orbits on a given level
or a sequence of levels converging to zero. The argument is particularly
transparent when $n=1$, i.e., for $M=\CP^1$.

\subsubsection{The construction for $T^*\CP^1$} Let $g$ be
the standard metric Hamiltonian on $W=T^*\CP^1\to \R$
equipped with the twisted symplectic structure. Note that all orbits
of the flow of $g$ on $W$ are periodic. Let $U$ be a neighborhood of
a level $g=a$. The universal cover $\tilde{U}$ of $U$ is 
symplectomorphic to a neighborhood of a sphere $S^3$ in $\C^2$, centered
at the origin. Let $\tilde{\Sigma}$ be 
an ellipsoid in $\tilde{U}$ which is close to $S^3$, invariant 
under deck transformations (the multiplication by $-1$), and carrying only 
finitely many closed characteristics. Then, $\tilde{\Sigma}$ descends to
a hypersurface $\Sigma$ in $U$ with a finite number of closed characteristics,
which is close to $g=a$. Note that $\tilde{\Sigma}$ can be taken
arbitrarily close to $S^3$ and hence $\Sigma$ can be made arbitrarily close
to $g=a$. Let now $\delta>0$ be small and $U=g^{-1}((a-\delta,a+\delta))$.
Then we can modify $g$ within $U$ so that the resulting function $G$ is
isotopic to $g$ but has $\Sigma$ as the level $G=a$. The next step is 
eliminating, as above, all periodic orbits on the level $G=a$. The resulting
function is now $C^2$-smooth.

Let us now pick a sequence $a_k\to 0$ and apply this process
to each $a=a_k$.  As a result, we obtain a continuous 
function $F$ which is $C^2$-smooth outside of the zero section and has
the same properties as the functions $F$ constructed in the two
previous proofs. As before, the proof for $n=1$ is concluded 
by applying Lemma \ref{lemma:smoothing}.

\begin{Remark}
A. Katok's example of a Finsler metric on $S^2$ with only two
closed geodesics can be easily described via a similar construction
of the hypersurface $\Sigma$, starting with the level $g=1$, in the standard
$T^*S^2$, and then symplectically embedding its double cover into $\C^2$.
\end{Remark}

\subsubsection{The general case: $T^*\CP^n$} Let $g_0$ be
the standard metric Hamiltonian on $W=T^*\CP^n$. As before, all orbits
of the flow of $g_0$ on $W$ are periodic and all orbits on a given
level have the same period.

Fix a level $g_0=a$. The first observation is that 
we can change $g_0$ in a neighborhood of $g_0=a$ to a function of the form 
$g=f\circ g_0$, without altering the level sets of the 
function or creating critical points, so that the new function $g$ is equal to 
$g_0$ outside the neighborhood and its flow has constant period (say, equal 
to one) near the level $g_0=a$. Without loss of generality, we may assume 
that $f(a)=a$, and hence the level we are working with is again $g=a$.

The next step is to notice that Ziller's method, \cite{zi}, applies to the 
function
$g$ near the level $g=a$, i.e., one can modify $g$ near the level to a 
$C^\infty$-function $G$, isotopic to $g$, so that the level $G=a$ carries
only a finite number of periodic orbits. 

For the sake of completeness, let us describe this modification in detail.
As in \cite{zi}, consider the action of $S^1=\R/\Z$ on $\CP^n$ 
induced by the following diagonal $S^1$-action on $\C^{n+1}$:
\begin{equation}
\labell{eq:S1-action}
t\cdot (z_0,\ldots,z_n)=(e^{2\pi\lambda_0 t}z_0,\ldots,e^{2\pi\lambda_n t}z_n),
\end{equation}
where $\lambda_0,\ldots,\lambda_n$ are mutually distinct integers and 
$\lambda_0=1$. This
action preserves $\sigma$ and $g_0$ and hence the lift $\psi_t$ of this
action to $T^*\CP^n$ preserves $g$ and the twisted symplectic structure
$\omega$. The flow $\psi_t$ is Hamiltonian.
Denote the Hamiltonian of $\psi_t$ by $g_1$ and set
$$
G=g+\alpha g_1,
$$
where $\alpha>0$ is a small irrational number. Note that since $\alpha$ is
small, the level $G=a$ lies near $g=a$ and, hence, is entirely contained in 
the region where the flow of $g$ has period one. We claim that the flow of 
$G$ on $T^*\CP^n$ has a finite number of periodic orbits on the level $G=a$
(and on nearby levels).

To prove this, let us first note that the closed orbits of the flow of $G$ 
near the 
level are in fact the closed orbits of $g$ which are invariant
under the flow $\psi_t$. This can be easily checked by using the
fact that the two flows commute and repeating word-for-word the reasoning 
from \cite[p. 138]{zi}. It follows that the projections of these periodic 
orbits of $g$ to $\CP^n$ are invariant under the $S^1$-action. 
Next observe that every orbit of $g$ is
a reparametrized orbit of the twisted geodesic flow on $g_0=a$. 

Hence, it 
suffices to show that the number of $S^1$-invariant twisted geodesics on 
$\CP^n$ coming from $g_0=a$ is finite. Let $\gamma$ be a twisted geodesic. 
The initial conditions $(\gamma(0),\gamma'(0))$ determine a projective line
$\CP^1\subset\CP^n$ and $\gamma$ is tangent to this $\CP^1$. Furthermore, 
this projective line is invariant under the $S^1$-action if $\gamma$ is 
invariant. (These facts can be seen as follows. It is well known that
$\gamma$ is the projection under the Hopf map of a great circle 
$\tilde{\gamma}$ in $S^{2n+1}$. The angle between $\tilde{\gamma}$ and the
Hopf fibers is
determined by $a$ and lies in the interval $(0,\, \pi/2)$. Thus
$\tilde{\gamma}(0)$ and $\tilde{\gamma}'(0)$ span a complex plane in 
$\C^{n+1}$. This complex plane gives rise to the projective line corresponding
to $\gamma$ and is invariant under the action \eqref{eq:S1-action} if
$\gamma$ is invariant.)

It is easy to
see that there are only finitely many (in fact, $n(n+1)/2$) projective 
lines invariant under the action \eqref{eq:S1-action}. (These 
projective lines arise from the $z_iz_j$-planes, $0\leq i< j\leq n$.) 
On such a projective line, the twisted geodesic $\gamma$ is an (oriented) 
$S^1$-invariant spherical circle whose geodesic curvature is determined by 
$a$. There are only two such oriented circles.

Thus, as we have shown, the flow of $G$ on the level $G=a$ has a finite
number of periodic orbits. Let us modify $G$ outside of a small neighborhood
of $G=a$ so that the resulting function, which we denote by $G$ again,
is $C^\infty$-smooth, equal to $g$ outside the neighborhood of $g=a$, and
isotopic to $g$.

The proof is finished in the same way as for $n=1$. Namely,
as the next step, we eliminate all periodic orbits on the level $G=a$. This
leads to a $C^\infty$-smooth  (since $\dim T^*\CP^n>4$) function which is
again isotopic to $g$. Next, we apply this construction to a sequence
$a_k\to 0$ of values of $g_0$ and obtain a continuous 
function $F$, $C^\infty$-smooth outside of the zero section
and it remains to again utilize Lemma \ref{lemma:smoothing}.
This concludes the proof of Proposition \ref{prop:non-exist-CPn} for all $n$.

\begin{Remark}
The above argument also shows, along the lines of \cite{zi}, that the flow 
of $G$ has exactly $n(n+1)$ periodic orbits on the level $G=a$.
\end{Remark}

\end{document}